\documentclass[11pt]{article}
\usepackage{fullpage} 
\usepackage[indent]{parskip} 

\usepackage[T1]{fontenc} 
\usepackage{lmodern} 
\usepackage{newtxtext}
\makeatletter%
\@ifpackageloaded{xeCJK}{}{%
    \usepackage[utf8]{inputenc}}
\makeatother

\usepackage[shortlabels,inline]{enumitem}

\usepackage{amsmath,amssymb,amsfonts}
\usepackage{dsfont}
\usepackage{bm}
\usepackage[mathscr]{eucal} 
\usepackage{amsthm}
\usepackage[mleftright]{diffcoeff}
	\difdef{f, s, fp, c}{}{%
		outer-Ldelim=\left.,%
		outer-Rdelim=\right|,%
		sub-nudge=0mu}%
	\difdef{c}{J}{%
		op-symbol=\rmJ,%
		op-sub-nudge=-2mu}%
\usepackage[group-separator={,}]{siunitx} 

\usepackage{xcolor}
\usepackage[pdftex]{graphicx} 
\usepackage{pgfplots}
	\pgfplotsset{compat=1.18}
	\usetikzlibrary{arrows}
\usepackage[export]{adjustbox}
\makeatletter%
\@ifclassloaded{IEEEtran}
	{\usepackage[caption=false,font=footnotesize]{subfig}}
	{\usepackage{subcaption}
	\captionsetup[figure]{font=footnotesize}}
\makeatother

\makeatletter%
\@ifpackageloaded{natbib}{}{%
    \usepackage[noadjust]{cite}} 
\makeatother

\usepackage{lipsum}
\usepackage{nameref}
\usepackage{silence}
\usepackage{verbatim} 
\usepackage{url} 
\usepackage{enotez}
	\setenotez{
		backref=true,
		counter-format=Roman,
		split=section,
		split-title={Notes for Section <ref>}
	}

\makeatletter
\newcommand{\disablepackage}[2]{%
  \disable@package@load{#1}{#2}}
\makeatother

    \setlist{nosep,leftmargin=*}

\theoremstyle{plain}
\newtheorem{thm}{Theorem}
\newtheorem{prop}{Proposition}
\newtheorem{cor}{Corollary}
\newtheorem{lem}{Lemma}

\theoremstyle{definition}

\newtheorem{eg}{Example}

\theoremstyle{remark}
\newtheorem{rmk}{Remark}

\newtheorem*{clm*}{Claim}

\makeatletter%
\@tempcnta=\@ne
\def\@nameedef#1{\expandafter\edef\csname #1\endcsname}
\loop\ifnum\@tempcnta<27
    \@nameedef{bb\@Alph\@tempcnta}{{\noexpand\mathbb{\@Alph\@tempcnta}}}
    \@nameedef{bf\@Alph\@tempcnta}{{\noexpand\mathbf{\@Alph\@tempcnta}}}
    \@nameedef{bm\@Alph\@tempcnta}{{\noexpand\bm{\@Alph\@tempcnta}}}
    \@nameedef{cal\@Alph\@tempcnta}{{\noexpand\mathcal{\@Alph\@tempcnta}}}
    \@nameedef{scr\@Alph\@tempcnta}{{\noexpand\mathscr{\@Alph\@tempcnta}}}
    \@nameedef{ds\@Alph\@tempcnta}{{\noexpand\mathds{\@Alph\@tempcnta}}}
    \@nameedef{frak\@Alph\@tempcnta}{{\noexpand\mathfrak{\@Alph\@tempcnta}}}
    \@nameedef{rm\@Alph\@tempcnta}{{\noexpand\mathrm{\@Alph\@tempcnta}}}
    \@nameedef{sf\@Alph\@tempcnta}{{\noexpand\mathsf{\@Alph\@tempcnta}}}
    \@nameedef{op\@Alph\@tempcnta}{{\noexpand\operatorname{\@Alph\@tempcnta}}}
    \@nameedef{hat\@Alph\@tempcnta}{{\noexpand\hat{\@Alph\@tempcnta}}}
    \@nameedef{tilde\@Alph\@tempcnta}{{\noexpand\tilde{\@Alph\@tempcnta}}}
    \@nameedef{ul\@Alph\@tempcnta}{{\noexpand\underline{\@Alph\@tempcnta}}}
    \@nameedef{ol\@Alph\@tempcnta}{{\noexpand\overline{\@Alph\@tempcnta}}}
	\@nameedef{v\@Alph\@tempcnta}{{\noexpand\vec{\@Alph\@tempcnta}}}
    \advance\@tempcnta\@ne
\repeat
\makeatother

\let\vec\undefined
\newcommand*{\vec}[1]{\bm{#1}}

\newcommand*{\Z}{\dsZ}

\newcommand*{\R}{\dsR}

\let\ForAll\forall
\let\forall\undefined
\DeclareMathOperator\forall{\ForAll}
\let\Exists\exists
\let\exists\undefined
\DeclareMathOperator\exists{\Exists}

\DeclareMathOperator{\tr}{tr}

\usepackage[colorlinks,citecolor=cyan]{hyperref}

\newcommand*{\oneMx}{\bm{1}}
\newcommand*{\pos}[1]{\max\left\{#1,\, 0\right\}}
\newcommand*{\card}[1]{|#1|}
\renewcommand*{\jacob}[2]{\difc.J.{#2}{#1}}
\newcommand*{\pdv}[2]{\difcp{#2}{#1}}
\newcommand*{\rpdv}[2]{\difcp[+]{#2}{#1}}
\newcommand*{\local}{}
\newcommand*{\Network}{\mathrm N}
\newcommand*{\Global}{\mathrm G}
\newcommand*{\lin}{\mathrm{lin}}
\newcommand*{\reach}{{}}
\newcommand*{\reachlim}{{}}
\newcommand*{\ini}{\mathrm{ini}}
\newcommand*{\iniinf}{*}

\DeclareMathOperator{\vol}{vol}
\DeclareMathOperator{\co}{co}
\DeclareMathOperator{\spabs}{\alpha}

\title{Topological entropy of nonlinear time-varying systems}
\author{Guosong Yang%
	\thanks{Guosong Yang is with the Department of Electrical and Computer Engineering, Rutgers University--New Brunswick, Piscataway, NJ 08854 USA (email: guosong.yang@rutgers.edu).}
\ and 
Daniel~Liberzon%
	\thanks{Daniel Liberzon is with the Coordinated Science Laboratory, University of Illinois Urbana-Champaign, Urbana, IL 61801 USA (e-mail: liberzon@illinois.edu).}}
\date{}

\begin{document}
\maketitle

\begin{abstract}
We study the topological entropy of nonlinear time-varying systems.
We establish conditions under which entropy is independent of the initial time and derive two general upper bounds based on the matrix measure of the system Jacobian and the spectral abscissa of a Metzler comparison matrix, respectively.
A general lower bound is also obtained using the trace of the Jacobian.
For interconnected systems, we first derive an upper bound using the matrix measure of an interconnection matrix constructed from the Jacobian blocks.
We then establish a new upper bound based on the spectral abscissa of a Metzler comparison matrix, incorporating the interconnection structure through subsystem-wise trajectory-separation analysis.
This bound also recovers the two general upper bounds as special cases.
All bounds depend on asymptotic properties of the Jacobian over reachable sets or their convex hulls.
Simple examples show that either general upper bound may be tighter, while a time-varying Lorenz system illustrates the effects of norm selection and subsystem decomposition.
\end{abstract}

\section{Introduction}\label{sec:intro}
In systems theory, topological entropy characterizes the rate at which information is generated in deterministic systems, captured by the increase in distinguishable behaviors at finite resolution, or by the complexity growth of a system acting on a set with finite measure.
Adler et al. \cite{AdlerKonheimMcAndrew1965} introduced the concept of topological entropy as an extension of Kolmogorov's metric entropy \cite{Kolmogorov1958}, quantifying the exponential growth rate of the minimal number of sets in iterated open-cover refinements.
An alternative definition, formulated in terms of the maximal number of separated trajectories over finite time horizons, was proposed by Dinaburg \cite{Dinaburg1970} and independently by Bowen \cite{Bowen1971}.
The equivalence between these two formulations was established in \cite{Bowen1971b}.
The variational principle, which identifies topological entropy as the supremum of metric entropy over invariant measures, was established through the combined works of Dinaburg \cite{Dinaburg1970,Dinaburg1971}, Goodman \cite{Goodman1971}, and Goodwyn \cite{Goodwyn1971}.
For a comprehensive introduction to topological entropy, see, e.g., \cite{KatokHasselblatt1995,Downarowicz2011} and the references therein.

Most existing results on topological entropy focus on time-invariant systems, while extending these results to time-varying dynamics requires accounting for the time dependence of trajectories and reachable sets.
Kolyada and Snoha \cite{KolyadaSnoha1996} defined topological entropy for discrete-time systems modeled by a sequence of selfmaps on a compact topological space, and established basic properties which strengthen under assumptions such as equicontinuity or uniform convergence.
Their definition and results were extended to time-varying maps defined on a sequence of compact metric spaces in \cite{KolyadaMisiurewiczSnoha1999}.
A variational principle based on this extended definition of topological entropy and a time-varying notion of metric entropy was established in \cite{KawanLatushkin2016}.
In the continuous-time setting, an upper bound on the topological entropy of a time-varying system on a compact invariant set was derived via the direct Lyapunov method in \cite{PogromskyMatveev2011}.

Topological entropy notions play a fundamental role in control theory, where feedback depends on the flow of information between sensors and actuators.
Nair et al. \cite{NairEvansMareelsMoran2004} introduced topological feedback entropy for discrete-time systems, building on the cover-refinement framework of \cite{AdlerKonheimMcAndrew1965}. Their definition extends classical topological entropy to quantify the exponential growth rate of control complexity required to ensure set invariance and stabilization.
For continuous-time systems, Colonius and Kawan \cite{ColoniusKawan2009} proposed invariance entropy, which is more closely aligned with the trajectory-counting formulation of \cite{Dinaburg1970,Bowen1971}.
The equivalence between topological feedback entropy and invariance entropy was established in \cite{ColoniusKawanNair2013}.
Suitable entropy notions have also been developed for exponential stabilization \cite{Colonius2012}, state estimation \cite{Savkin2006,MatveevPogromsky2016,LiberzonMitra2018}, and model detection \cite{LiberzonMitra2018}.

Entropy-based notions are closely linked to data-rate requirements for control.
For continuous-time linear time-invariant systems, the minimal data rate for feedback stabilization is given by the sum of the positive real parts of the eigenvalues of the system matrix, counted with algebraic multiplicity \cite{HespanhaOrtegaVasudevan2002}; 
in the discrete-time case, it is the sum of their logarithms \cite{HespanhaOrtegaVasudevan2002,NairEvans2003,TatikondaMitter2004}.
This quantity coincides with the topological entropy of the open-loop system \cite{Bowen1971}.
While this connection has also been extensively studied for nonlinear time-invariant systems (e.g., \cite{NairEvansMareelsMoran2004,LiberzonHespanha2005,ColoniusKawan2009,Colonius2012}), its counterpart for time-varying dynamics remains largely unexplored, except in certain cases such as switched systems, discussed below.

In systems and control, a common paradigm is to decompose a complex system into an interconnection of simpler subsystems and derive system-level properties from those of the subsystems.
A classical example is the small-gain theorem, used to establish stability of interconnected linear and nonlinear systems (e.g., \cite{DesoerVidyasagar2009,LiuJiangHill2014}).
Entropy-based approaches have also been applied in this context, particularly for discrete-time systems.
The invariance entropy of a network of control systems was related to the entropy of its subsystems and associated data rates in \cite{KawanDelvenne2016}.
Observability rates for networked systems were linked to topological entropy and auxiliary constructs such as storage functions and supply rates derived from subsystem linearization in \cite{MatveevProskurnikovPogromskyFridman2019}.
Invariance feedback entropy for a network of uncertain systems was investigated in \cite{TomarZamani2020}.
In the continuous-time setting, an upper bound on the topological entropy of interconnected systems was derived in \cite{Liberzon2021}, which can be viewed as a prototype for some of the results in the present paper.
However, much remains to be understood about when and how explicit bounds on the topological entropy of an interconnected nonlinear time-varying system can be deduced from those of its subsystems.

Topological entropy and data-rate requirements for control have been studied for switched systems, an important class of time-varying systems.
Sufficient data rates for feedback stabilization of switched linear systems were established in \cite{Liberzon2014,YangLiberzon2018}.
Similar data-rate conditions were derived in \cite{SibaiMitra2017} by extending a notion of estimation entropy from \cite{LiberzonMitra2018}.
Formulae and bounds on the topological entropy of continuous-time switched linear systems were derived in \cite{YangSchmidtLiberzon2018,YangHespanha2018,YangSchmidtLiberzonHespanha2020}, with their relationship to stability conditions analyzed in \cite{YangHespanhaLiberzon2019,YangSchmidtLiberzonHespanha2020}.
In an extension to switched nonlinear and interconnected systems \cite{YangLiberzonHespanha2023}, entropy bounds were constructed using the matrix measure or the trace of the system Jacobian. These quantities also appear in the bounds developed in the present paper.
For discrete-time switched linear systems, the topological entropy under worst-case switching sequences was studied using the joint spectral radius in \cite{BergerJungers2020}, and the connection between estimation entropy, Lyapunov exponents, and average data rate was examined under suitable regularity conditions in \cite{VicinansaLiberzon2023}.

The main objective of this paper is to establish upper and lower bounds on the topological entropy of general and interconnected nonlinear time-varying systems. 

We introduce the definition of topological entropy in Section~\ref{sec:pre} and examine its dependence on the initial time and the set of initial states. 
We show that changing the initial time does not affect the entropy if the system is reinitialized from the corresponding reachable set, but it may differ if initialized from the original initial set.

In Section~\ref{sec:tv-ent}, we establish two general upper bounds on the entropy of nonlinear time-varying systems. 
One upper bound is constructed using the matrix measure of the system Jacobian, while the other is based on the spectral abscissa of a Metzler comparison matrix.
A general lower bound is also established using the trace of the Jacobian.
We describe the use of possibly time-varying enclosures of reachable sets to obtain more tractable bounds, potentially at the expense of tightness.
Simple examples illustrate that, for a fixed induced norm, either upper bound may be tighter.

In Section~\ref{sec:inter-ent}, we consider interconnected nonlinear time-varying systems.
We first adapt the general upper bound based on matrix measure from Section~\ref{sec:tv-ent}, using an interconnection matrix constructed from the Jacobian blocks. 
We then develop a new upper bound based on the spectral abscissa of a Metzler comparison matrix associated with the interconnection. 
This bound incorporates the interconnection structure through subsystem-wise trajectory-separation analysis.
It recovers both upper bounds from Section~\ref{sec:tv-ent} by viewing the general system as an interconnection composed of either a single subsystem or scalar subsystems.
All entropy bounds depend on asymptotic properties of the Jacobian over reachable sets or their convex hulls.

In Section~\ref{sec:alt}, we provide alternative upper bounds on entropy for both general and interconnected systems.
Section~\ref{sec:eg} illustrates the effects of norm selection and subsystem decomposition on the entropy bounds using a time-varying Lorenz system.
The proofs of the main results are presented in Section~\ref{sec:proof}, with some technical arguments deferred to the appendices. 
Section~\ref{sec:end} concludes the paper with a brief summary and remarks on future research directions.

\emph{Notations:}
For vectors $ v_1 \in \R^{n_1}, \ldots, v_m \in \R^{n_m} $, let $ (v_1, \ldots, v_m) \in \R^{n_1 + \cdots + n_m} $ denote their concatenation.
For a square matrix $ M \in \R^{n \times n} $, let $ \tr(M) $ and $ \spabs(M) $ denote its trace and spectral abscissa (the largest real part of its eigenvalues), respectively.
Let $ \card{E} $ denote the cardinality of a finite set $ E $.
For a set $ S \subset \R^n $, let $ \co(S) $ and $ \vol(S) $ denote its convex hull and volume (Lebesgue measure), respectively.
By default, all logarithms are natural logarithms (to avoid an extra multiplicative factor $ \ln 2 $ in entropy computations).

\section{Preliminaries}\label{sec:pre}
Consider the continuous-time \emph{nonlinear time-varying} dynamical system
\begin{equation}\label{eq:tv}
    \dot x = f(t, x), \qquad t \geq t_0,
\end{equation}
where $ x \in \R^n $ is the state and $ t_0 \in \R $ is the initial time. 
Assume that the system function $ f(t, x) $ is continuously differentiable in $ x $, 
that $ f(t, x) $ and the system Jacobian $ \jacob{x}{f(t, x)} $ are both piecewise continuous in $ t $,%
\footnote{Here, piecewise continuity in $ t $ is understood locally uniformly in $ x $, with a locally finite set of discontinuity times independent of $ x $, and finite one-sided limits at each discontinuity.}
and that \eqref{eq:tv} is forward complete.
Let $ \xi(t, t_0, x) $ denote the solution of \eqref{eq:tv} at time $ t $ with initial state $ x $ at time $ t_0 $. 
Under these assumptions, $ \xi(t, t_0, x) $ is unique, continuously differentiable in $ x $, absolutely continuous in $ t $, and satisfies \eqref{eq:tv} away from discontinuities of $ f(t, x) $ in $ t $. 
Moreover, the solution Jacobian $ \jacob{x}{\xi(t, t_0, x)} $ is absolutely continuous and piecewise continuously differentiable in $ t $. 
For brevity, we use the notation
\[
	\jacob{x}{f(t, v)}
	:= \jacob{x}{f(t, x)}[x=v], 
	\quad 
	\jacob{x}{\xi(t, t_0, v)} 
	:= \jacob{x}{\xi(t, t_0, x)}[x=v]
\]
for the system and solution Jacobians evaluated at $ x = v $, respectively.

\subsection{Entropy definition}\label{ssec:pre-ent}
We define the topological entropy of the nonlinear time-varying system \eqref{eq:tv} with initial states in a compact set $ K \subset \R^n $ with nonempty interior, called the \emph{initial set}. 
Let $ |\cdot| $ be a norm on $ \R^n $, and $ \|\cdot\| $ the corresponding induced norm on $ \R^{n \times n} $. 
Given a time horizon $ T \geq 0 $ and a radius $ \varepsilon > 0 $, define the following open ball in $ \R^n $ centered at $ x $:
\begin{align}
     &B_{f, t_0}(x, \varepsilon, T) \notag\\
     &:= \left\{ \bar x \in \R^n: \max_{t \in [t_0, t_0 + T]} |\xi(t, t_0, \bar x) - \xi(t, t_0, x)| < \varepsilon \right\}. \label{eq:ball-dfn}
\end{align}
A subset $ E \subset K $ is said to be \emph{$ (T, \varepsilon) $-spanning} (for the initial set $ K $ at initial time $ t_0 $) if
\begin{equation}\label{eq:span-dfn}
    K \subset \bigcup_{x \in E} B_{f, t_0}(x, \varepsilon, T),
\end{equation}
that is, for every $ \bar x \in K $, there exists $ x \in E $ such that
\[
	|\xi(t, t_0, \bar x) - \xi(t, t_0, x)| < \varepsilon 
	\qquad
	\forall t \in [t_0, t_0 + T]. 
\]
A $ (T, \varepsilon) $-spanning set of minimum cardinality is called minimal.
Its cardinality, denoted by $ S(f, t_0, \varepsilon, T, K) \geq 1 $, is nondecreasing in $ T $ and nonincreasing in $ \varepsilon $.
The \emph{topological entropy} of \eqref{eq:tv} with initial set $ K $ at initial time $ t_0 $ is defined in terms of the exponential growth rate of $ S(f, t_0, \varepsilon, T, K) $:
\begin{equation}\label{eq:ent-dfn}
    h(f, t_0, K) := \lim_{\varepsilon \searrow 0} \limsup_{T \to \infty} \frac{1}{T} \log S(f, t_0, \varepsilon, T, K).
\end{equation}
For brevity, we occasionally refer to $ h(f, t_0, K) $ simply as the (topological) entropy of \eqref{eq:tv}.
By the equivalence of norms on finite-dimensional vector spaces, the value of $ h(f, t_0, K) $ is independent of the choice of norm $ |\cdot| $ on $ \R^n $. 
In particular, it is invariant under a change of basis.
 
An equivalent definition of entropy can be formulated using separated sets; see \cite{YangLiberzonHespanha2023} for its application to switched nonlinear systems.

When \eqref{eq:tv} is time-invariant, we adopt the convention $ t_0 = 0 $ and consider the \emph{nonlinear time-invariant} system
\begin{equation}\label{eq:ti}
	\dot x = f(x), \qquad t \geq 0,
\end{equation}
where $ f: \R^n \to \R^n $ is a continuously differentiable function.
For time-invariant systems, we omit $ t_0 $ from the entropy-related notation introduced above.

\begin{rmk}\label{rmk:liminf}
For time-invariant systems with initial sets contained in compact forward-invariant sets, the upper limit in \eqref{eq:ent-dfn} can be replaced by the lower limit \cite[pp.~109--110]{KatokHasselblatt1995}. 
For time-varying systems, the subadditivity underlying this equivalence (see the proof of \cite[Lemma~3.1.5, p.~109]{KatokHasselblatt1995}) may fail, so the same argument does not directly apply.
\end{rmk}

\subsection{Initial set and initial time}\label{ssec:pre-ini}
We examine how the topological entropy $ h(f, t_0, K) $ of the nonlinear time-varying system \eqref{eq:tv} depends on the initial set $ K $ and initial time $ t_0 $.

Let $ \xi(t, t_0, K) := \{\xi(t, t_0, x): x \in K\} $ denote the \emph{reachable set} at time $ t $ starting from $ K $ at time $ t_0 $. In particular, $ \xi(t_0, t_0, K) = K $. 
The following lemma shows that advancing the initial time from $ t_0 $ to $ t_1 $, while replacing $ K $ with the reachable set at $ t_1 $, does not change the entropy. 
The proof is given in Appendix~\ref{apx:ini-time}.
\begin{lem}\label{lem:ini-time}
The topological entropy of \eqref{eq:tv} satisfies
\begin{equation}\label{eq:ini-time}
	h(f, t_0, K) = h(f, t_1, \xi(t_1, t_0, K)) \qquad \forall t_1 \geq t_0.
\end{equation}
\end{lem}

Related results for discrete-time systems with equicontinuous functions and for continuous-time systems on compact invariant sets were established in \cite[Prop.~2.1]{KawanLatushkin2016} and \cite[Lemma~2.8]{PogromskyMatveev2011}, respectively, building on the time-varying entropy framework of \cite{KolyadaSnoha1996,KolyadaMisiurewiczSnoha1999}.

However, the entropy of \eqref{eq:tv} may depend on the initial set and, even for a fixed initial set, on the initial time, as illustrated by the following example.
\begin{eg}\label{eg:ini-time-diff}
Consider the continuously differentiable function $ g: \R \to \R $ defined by
\begin{equation*}
	g(x) := \begin{cases}
		0, &x < -1, \\
		\sqrt{3} - \sqrt{3 - (x+1)^2}, &-1 \leq x \leq 1/2, \\
		\sqrt{3} x, &x > 1/2.
	\end{cases}
\end{equation*}
\begin{enumerate}
	\item (Different initial sets). 
		Consider the nonlinear time-invariant system \eqref{eq:ti} with $ f(x) = g(x) $. 
		For initial sets $ K_1 \subset (-\infty, -1) $ and $ K_2 \subset (1, \infty) $, we have $ h(g, K_1) = 0 $ and $ h(g, K_2) = \sqrt{3} $.
	\item (Same initial set at different initial times).
		Consider the nonlinear time-varying system \eqref{eq:tv} with
		\begin{equation*}
			f(t, x) := \begin{cases}
				g(x+4), &t < 0, \\
				g(x), &t \geq 0,
			\end{cases}
		\end{equation*}
		and fix the initial set $ K := [-3, -2] $. 
		For an initial time $ t_0 \geq 0 $, we have $ h(f, t_0, K) = h(g, K) = 0 $. 
		For an initial time $ t_1 \leq -1 $ and $ x \in K $, the solution satisfies 
		$ \xi(t, t_1, x) = e^{\sqrt{3} (t - t_1)} (x + 4) - 4 $
		for $ t \in [t_1, 0) $. 
		Thus, the reachable set at time $ 0 $ is
		$ \xi(0, t_1, K) = [e^{-\sqrt{3} t_1} - 4,\, 2 e^{-\sqrt{3} t_1} - 4] \subset (1, \infty) $. 
		By \eqref{eq:ini-time}, we have 
		$ h(f, t_1, K) = h(f, 0, \xi(0, t_1, K)) = h(g, \xi(0, t_1, K)) = \sqrt{3} $.
\end{enumerate}
\end{eg}

Consider the special case of \eqref{eq:tv} given by the \emph{linear time-varying (LTV)} system
\begin{equation}\label{eq:ltv}
    \dot x = A(t)\, x, \qquad t \geq t_0,
\end{equation}
where $ A: \R \to \R^{n \times n} $ is a piecewise-continuous matrix-valued function.
The following lemma extends \cite[Prop.~2]{YangSchmidtLiberzonHespanha2020} from switched linear systems to general LTV systems.
\begin{lem}\label{lem:ini-ltv}
The topological entropy of the LTV system \eqref{eq:ltv} is independent of the choices of initial set and initial time.
\end{lem}

\begin{proof}
Independence from the initial set follows by an argument analogous to the proof of \cite[Prop.~2]{YangSchmidtLiberzonHespanha2020}, which relies only on linearity of the solution map and therefore applies to general LTV systems. 
Independence from the initial time then follows from Lemma~\ref{lem:ini-time}.
\end{proof}

In light of Lemma~\ref{lem:ini-ltv}, we let $ h(A) $ denote the entropy of the LTV system \eqref{eq:ltv}.

\section{Entropy of nonlinear time-varying systems}\label{sec:tv-ent}
In this section, we present upper and lower bounds on the topological entropy of the nonlinear time-varying system \eqref{eq:tv}. 

We first recall the notion of matrix measure, which will be used to derive upper bounds on entropy.
Given an induced norm $ \|\cdot\| $ on $ \R^{n \times n} $, the \emph{matrix measure} (or \emph{logarithmic norm}) of $ A \in \R^{n \times n} $ is defined by
\[
    \mu(A) := \lim_{t \searrow 0} \frac{\|I + t A\| - 1}{t}.
\]
As shown in \cite[Th.~2.2.16, p.~22]{Vidyasagar2002}, the matrix measure satisfies 
\begin{equation}\label{eq:meas-ineq-eig}
    \spabs(A) 
    \leq \mu(A) 
    \leq \|A\| 
    \qquad 
    \forall A \in \R^{n \times n}.
\end{equation}
Moreover,
\begin{equation*}\label{eq:meas-ineq-sum}
	\mu(A + B) \leq \mu(A) + \mu(B), 
	\qquad 
	\mu(cA) = c \mu(A) 
\end{equation*}
for all $ A, B \in \R^{n \times n} $ and $ c \geq 0 $.
In particular, $ \mu(\cdot) $ is convex. For standard induced norms, explicit expressions for the matrix measure are available. 
For example, the matrix measures of $ A = [a_{ij}] \in \R^{n \times n} $ associated with the induced $ \infty $- and $ 2 $-norms are
\begin{equation*}\label{eq:meas-inf-norm}
    \mu_\infty(A) 
    = \max_{1 \leq i \leq n} \Bigg( 
    	a_{ii} + \sum_{j \neq i} |a_{ij}| 
    \Bigg),
    \quad
    \mu_2(A) = \alpha \left( 
    	\frac{A + A^\top}{2}
    \right),
\end{equation*}
respectively. 
See \cite[Sec.~2.2.2]{Vidyasagar2002} and \cite[Fact~11.15.7, p.~690]{Bernstein2009} for further properties of the matrix measure.

The following theorem establishes upper and lower bounds on the entropy of \eqref{eq:tv}. 
The proof, which partly relies on results for interconnected systems developed in the next section, is given in Section~\ref{ssec:tv-ent-bnd-pf}.
\begin{thm}\label{thm:tv-ent-bnd}
The topological entropy of the nonlinear time-varying system \eqref{eq:tv} is upper bounded by
\begin{equation}\label{eq:tv-ent-upper-reach}
    h(f, t_0, K) \leq n \pos{\hat\mu^\reachlim},
\end{equation}
where
\begin{equation}\label{eq:tv-meas-lim-reach}
	\hat\mu^\reachlim 
	:= \limsup_{t \to \infty} 
	\max_{v \in \co(\xi(t, t_0, K))} 
	\mu(\jacob{x}{f(t, v)}),
\end{equation}
and lower bounded by
\begin{equation}\label{eq:tv-ent-lower}
    h(f, t_0, K) \geq \pos{\check\chi},
\end{equation}
where
\begin{equation}\label{eq:tv-tr-lim}
	\check\chi 
	:= \liminf_{t \to \infty} 
	\min_{v \in \xi(t, t_0, K)} 
	\tr(\jacob{x}{f(t, v)}).
\end{equation}
\end{thm}

The upper bound in \eqref{eq:tv-ent-upper-reach} holds for any induced norm used to compute the matrix measure in \eqref{eq:tv-meas-lim-reach}, and may thus be improved by a suitable choice of norm.

The upper bound in \eqref{eq:tv-ent-upper-reach} depends on the asymptotic behavior of the matrix measure of the Jacobian $ \jacob{x}{f(t, x)} $ over the convex hull of the reachable set $ \xi(t, t_0, K) $.
Similarly, the lower bound in \eqref{eq:tv-ent-lower} depends on the asymptotic behavior of the trace of the Jacobian  over the reachable set.
These formulations allow the bounds to be evaluated using suitable enclosures of the reachable sets, as discussed in the following remark.

\begin{rmk}\label{rmk:superset}
The entropy bounds in Theorem~\ref{thm:tv-ent-bnd} can be evaluated using more tractable enclosures of the reachable sets, potentially at the expense of tightness.
For example, consider a set $ S $ satisfying any of the following conditions:
\begin{itemize}
	\item $ S = \R^n $;
	\item $ \co(\xi(t, t_0, K)) \subset S $ for all sufficiently large $ t $; or
	\item $ S $ is convex, forward invariant, and satisfies $ \co(\xi(t_1, t_0, K)) \subset S $ for some $ t_1 \geq t_0 $ (e.g., $ \co(K) \subset S $).
\end{itemize}
Then the entropy bounds can be simplified by bounding the relevant quantities over $ S $. 
For instance, if a constant $ \hat\mu^S $ satisfies $ \hat\mu^S \geq \mu(\jacob{x}{f(t, v)}) $ for all $ v \in S $ and sufficiently large $ t $, then the upper bound in \eqref{eq:tv-ent-upper-reach} remains valid with $ \hat\mu^S $ replacing $ \hat\mu^\reachlim $. 
More generally, one may consider a forward-invariant family of convex sets $ S_t $ satisfying $ \co(\xi(t_1, t_0, K)) \subset S_{t_1} $ for some $ t_1 \geq t_0 $, and bound the relevant quantities over $ S_t $.
The use of such time-varying enclosures is illustrated in Example~\ref{eg:tv-lorenz} in Section~\ref{sec:eg}.
Analogous substitutions apply to other entropy bounds presented in this paper.
\end{rmk}

\begin{rmk}\label{rmk:stability}
The upper bound in Theorem~\ref{thm:tv-ent-bnd} implies that the topological entropy of \eqref{eq:tv} is zero whenever $ \hat\mu^\reachlim \leq 0 $.
This includes systems that are contracting with respect to the same norm, for which trajectories converge toward one another exponentially \cite{Sontag2010,RussoBernardoSontag2013}.
More generally, asymptotic stability can also imply zero topological entropy for initial sets contained in the basin of attraction.
Conversely, zero topological entropy does not imply that the upper bounds derived in this paper vanish or that the system is stable.
For the \emph{linear time-invariant (LTI)} system
\begin{equation}\label{eq:lti}
	\dot x = A x,
\end{equation}
the topological entropy is known to equal the sum of the positive real parts of the eigenvalues of $ A $, counted with algebraic multiplicity (see, e.g., \cite[Lemma~3]{YangSchmidtLiberzonHespanha2020}).
For
\[
	A = \begin{bmatrix}
		1 & 1 \\
		-1 & -1
	\end{bmatrix},
\]
since $ A^2 = 0 $ and $ e^{At} = I + tA $, the LTI system \eqref{eq:lti} has zero entropy but is not stable.
See \cite[Sec.~7]{YangSchmidtLiberzonHespanha2020} and \cite{YangHespanhaLiberzon2019} for further discussions on entropy and stability of switched linear systems.
\end{rmk}

As a preview of the results for interconnected systems, we derive an alternative upper bound on the entropy of \eqref{eq:tv} by viewing it as $ n $ interconnected scalar subsystems. 
Unlike \eqref{eq:tv-ent-upper-reach}, this bound involves neither a matrix measure nor an induced norm.
The proof is given in Section~\ref{ssec:tv-ent-bnd-pf}. 
For each $ i \in \{1, \ldots, n\} $, let $ x_i $ and $ f_i(t, x) $ denote the $ i $-th components of the state $ x $ and the vector field $ f(t, x) $ in \eqref{eq:tv}, respectively.
\begin{prop}\label{prop:tv-ent-upper-eig-reach}
Let $ \hat A^\reachlim = [\hat a^\reachlim_{ij}] \in \R^{n \times n} $ be any matrix satisfying
\begin{equation}\label{eq:tv-comp-lim-reach}
\begin{aligned}
	\hat a^\reachlim_{ii} 
	&\geq \limsup_{t \to \infty} 
	\max_{v \in \co(\xi(t, t_0, K))} 
	\pdv{x_i}{f_i(t, v)},
	\\
	\hat a^\reachlim_{ij} 
	&\geq \limsup_{t \to \infty} 
	\max_{v \in \co(\xi(t, t_0, K))} 
	|\pdv{x_j}{f_i(t, v)}|
\end{aligned}
\end{equation}
for all diagonal and off-diagonal entries, respectively.
Then the topological entropy of the nonlinear time-varying system \eqref{eq:tv} is upper bounded in terms of the spectral abscissa of $ \hat A^\reachlim $:
\begin{equation}\label{eq:tv-ent-upper-eig-reach}
    h(f, t_0, K) \leq n \pos{\spabs(\hat A^\reachlim)}.
\end{equation}
\end{prop}

Note that any matrix satisfying \eqref{eq:tv-comp-lim-reach} is Metzler (i.e., has nonnegative off-diagonal entries).

The two upper bounds in \eqref{eq:tv-ent-upper-reach} and \eqref{eq:tv-ent-upper-eig-reach} are obtained by applying different operations to the Jacobian $ \jacob{x}{f(t, x)} $.
The former evaluates a matrix measure before taking the asymptotic maximum, 
whereas the latter first constructs a constant Metzler comparison matrix from asymptotic entrywise bounds and then computes its spectral abscissa.
For a fixed choice of induced norm, either bound may be tighter, even for linear systems, as illustrated by the following examples using the induced $ \infty $-norm.
The effect of norm selection is further illustrated using a time-varying Lorenz system in Example~\ref{eg:tv-lorenz} in Section~\ref{sec:eg}.

\begin{eg}\label{eg:ltv-ent-upper}
Consider the LTV system \eqref{eq:ltv} with 
\[ 
	A(t) = \begin{bmatrix}
		\sin(t) & \cos(t) \\
		\sin(t) & \cos(t)
	\end{bmatrix}. 
\]
For the induced $ \infty $-norm, 
\[ 
	\mu_\infty(A(t)) = \max\{\sin(t) + |{\cos(t)}|,\, |{\sin(t)}| + \cos(t)\},
\]
and hence $ \limsup_{t \to \infty} \mu_\infty(A(t)) = \sqrt{2} $. 
Thus, \eqref{eq:tv-ent-upper-reach} yields $ h(A) \leq 2 \sqrt{2} $.
On the other hand, the smallest Metzler matrix satisfying \eqref{eq:tv-comp-lim-reach} is
\[
	\hat A = \begin{bmatrix}
		1 & 1 \\
		1 & 1
	\end{bmatrix},
\]
for which $ \spabs(\hat A) = 2 $.
Thus, \eqref{eq:tv-ent-upper-eig-reach} yields $ h(A) \leq 4 $.
Therefore, the upper bound in \eqref{eq:tv-ent-upper-reach} is tighter for this example.
\end{eg}

\begin{eg}\label{eg:lti-ent-upper}
Consider the LTI system \eqref{eq:lti} with the Metzler matrix
\[ 
	A = \begin{bmatrix}
		0 & 2 \\
		1 & 0
	\end{bmatrix}. 
\]
For the induced $ \infty $-norm, $ \mu_\infty(A) = 2 $, while $ \spabs(A) = \sqrt{2} $.
Therefore, \eqref{eq:tv-ent-upper-reach} yields $ h(A) \leq 4 $, 
whereas \eqref{eq:tv-ent-upper-eig-reach} with $ \hat A = A $ yields the tighter bound $ h(A) \leq 2 \sqrt{2} $ for this example.
For comparison, the exact entropy is $ h(A) = \sqrt{2} $.
\end{eg}

Although the upper bound in \eqref{eq:tv-ent-upper-reach} may be tightened by optimizing over induced norms, such optimization can be challenging for nonlinear time-varying systems.
The bound in \eqref{eq:tv-ent-upper-eig-reach} provides an alternative based on asymptotic entrywise bounds on the Jacobian.
Additional upper bounds on the entropy of \eqref{eq:tv} are given in Section~\ref{sec:alt}.

\section{Entropy of interconnected systems}\label{sec:inter-ent}
In this section, we consider the case where the nonlinear time-varying system \eqref{eq:tv} is composed of $ m $ interconnected subsystems. 
For each $ i \in \{1, \ldots, m\} $, let $ x_i \in \R^{n_i} $ denote the state of the $ i $-th subsystem, and $ f_i(t, x) \in \R^{n_i} $ the corresponding components of $ f(t, x) $. 
Then $ n_1 + \cdots + n_m = n $, and \eqref{eq:tv} is given by the \emph{interconnected} nonlinear time-varying system
\begin{equation}\label{eq:inter}
	\dot x_i = f_i(t, x_1, \ldots, x_m), 
	\quad 
	t \geq t_0, 
	\qquad 
	i \in \{1, \ldots, m\}.
\end{equation}
The overall state and vector field are then given by $ x = (x_1, \ldots, x_m) $ and $ f(t, x) = (f_1(t, x), \ldots, f_m(t, x)) $, respectively. 
For brevity, we use the notation
\[
	\jacob{x_j}{f_i(t, v)} := \jacob{x_j}{f_i(t, x)}[x=v]
\]
for the Jacobian block of $ f_i(t, x) $ with respect to $ x_j $ evaluated at $ x = v $. 

In general, the topological entropy of \eqref{eq:inter} cannot be determined solely from the entropies of its uncoupled subsystems. 
For example, the LTI system in Remark~\ref{rmk:stability} has zero entropy, while its uncoupled scalar subsystems have entropies $ 1 $ and $ 0 $. 
Conversely, the LTI system in Example~\ref{eg:lti-ent-upper} has entropy $ h(A) = \sqrt{2} $, while both uncoupled scalar subsystems have zero entropy.

We first derive an upper bound on the topological entropy of the interconnected system \eqref{eq:inter} by adapting the upper bound in Theorem~\ref{thm:tv-ent-bnd}, using matrix measures and induced norms of the Jacobian blocks. 

Following \cite{RussoBernardoSontag2013}, suppose the following norms are given:
\begin{enumerate}
	\item for each $ i \in \{1, \ldots, m\} $, a ``local'' norm $ |\cdot|_{\local i} $ on $ \R^{n_i} $, and
	\item a ``network'' norm $ |\cdot|_\Network $ on $ \R^m $ that is \emph{monotone}: for all nonnegative vectors $ v, w \in \R_{\geq 0}^m $,
		\begin{equation}\label{eq:mono}
			v \geq w \implies |v|_\Network \geq |w|_\Network.
		\end{equation}
		In particular, all $ p $-norms with $ p \geq 1 $ are monotone.%
		\footnote{Inequalities between vectors or matrices are interpreted entrywise (e.g., $ A \geq 0 $ means that $ A $ is a nonnegative matrix).}
\end{enumerate}
For a vector $ v = (v_1, \ldots, v_m) \in \R^n $ with $ v_i \in \R^{n_i} $, the ``global'' norm $ |\cdot|_\Global $ is defined by
\begin{equation}\label{eq:global-norm-dfn}
	|v|_\Global := |(|v_1|_{\local 1}, \ldots, |v_m|_{\local m})|_\Network.
\end{equation}
As $ |\cdot|_\Network $ is monotone, one can verify that $ |\cdot|_\Global $ satisfies the triangle inequality and is indeed a norm. 
Let $ \|\cdot\|_{\local i} $, $ \|\cdot\|_\Network $, and $ \|\cdot\|_\Global $ denote the corresponding induced norms on $ \R^{n_i \times n_i} $, $ \R^{m \times m} $, and $ \R^{n \times n} $, respectively, and $ \mu_{\local i}(\cdot) $, $ \mu_\Network(\cdot) $, and $ \mu_\Global(\cdot) $ the corresponding matrix measures. 
The corresponding induced norms $ \|\cdot\|_{\local ij} $ on $ \R^{n_i \times n_j} $ are defined by
\begin{equation*}
	\|A\|_{\local ij} := \max_{|v|_{\local j} = 1} |A v|_{\local i}, \qquad i, j \in \{1, \ldots, m\}.
\end{equation*}

We define the \emph{interconnection matrix} of \eqref{eq:inter} as the matrix-valued function $ A_\Network(t, x) = [a_{ij}(t, x)] \in \R^{m \times m} $ with
\begin{equation}\label{eq:inter-sub-meas}
	a_{ii}(t, x) 
	:= \mu_{\local i}(\jacob{x_i}{f_i(t, x)}), 
	\qquad 
	a_{ij}(t, x) 
	:= \|\jacob{x_j}{f_i(t, x)}\|_{\local ij}
\end{equation}
for diagonal and off-diagonal entries, respectively.
The function $ A_\Network(t, x) $ is referred to as a \emph{structure Jacobian} in \cite{RussoBernardoSontag2013}, where it is used to establish contraction conditions.

The following result uses the interconnection matrix $ A_\Network(t, x) $ to derive an upper bound on the entropy of \eqref{eq:inter} without explicitly evaluating a matrix measure of the full Jacobian $ \jacob{x}{f(t, x)} $.
\begin{cor}\label{cor:inter-ent-upper-reach}
Consider the interconnection matrix $ A_\Network(t, x) $ defined in \eqref{eq:inter-sub-meas}. 
The topological entropy of the interconnected system \eqref{eq:inter} is upper bounded by
\begin{equation}\label{eq:inter-ent-upper-reach}
    h(f, t_0, K) \leq n \pos{\hat\mu^\reachlim_\Network},
\end{equation}
where
\begin{equation}\label{eq:inter-meas-lim-reach}
	\hat\mu^\reachlim_\Network 
	:= \limsup_{t \to \infty} 
	\max_{v \in \co(\xi(t, t_0, K))} 
	\mu_\Network(A_\Network(t, v)).
\end{equation}
\end{cor}

\begin{proof}
By Lemma~\ref{lem:inter-meas-bnd} below, the constant $ \hat\mu^\reachlim $ defined in \eqref{eq:tv-meas-lim-reach} using the global matrix measure $ \mu_\Global(\cdot) $ satisfies $ \hat\mu^\reachlim \leq \hat\mu^\reachlim_\Network $. 
Thus, the upper bound in \eqref{eq:inter-ent-upper-reach} follows from \eqref{eq:tv-ent-upper-reach}.
\end{proof}

\begin{lem}[{\cite[Th.~2]{RussoBernardoSontag2013}}]\label{lem:inter-meas-bnd}
For a block matrix $ A = [A_{ij}] \in \R^{n \times n} $ with $ A_{ij} \in \R^{n_i \times n_j} $ for $ i, j \in \{1, \ldots, m\} $, let $ A_\Network = [a_{ij}] \in \R^{m \times m} $ be any matrix satisfying
\begin{equation*}
	a_{ii} \geq \mu_{\local i}(A_{ii}), 
	\qquad 
	a_{ij} \geq \|A_{ij}\|_{\local ij}
\end{equation*}
for all diagonal and off-diagonal entries, respectively.
Then
\begin{equation*}
	\mu_\Global(A) \leq \mu_\Network(A_\Network).
\end{equation*}
\end{lem}

We now present a new upper bound on the entropy of \eqref{eq:inter}. 
The proof is given in Section~\ref{ssec:inter-ent-upper-eig-reach-pf}, using preliminary results established in Sections~\ref{ssec:soln-vol} and~\ref{ssec:ent-span}.
\begin{thm}\label{thm:inter-ent-upper-eig-reach}
Consider the interconnection matrix $ A_\Network(t, x) = [a_{ij}(t, x)] \in \R^{m \times m} $ defined in \eqref{eq:inter-sub-meas}. 
Let $ \hat A^\reachlim_\Network = [\hat a^\reachlim_{ij}] \in \R^{m \times m} $ be any matrix satisfying
\begin{equation}\label{eq:inter-sub-meas-lim-reach}
	\hat a^\reachlim_{ij} 
	\geq \limsup_{t \to \infty} 
	\max_{v \in \co(\xi(t, t_0, K))} 
	a_{ij}(t, v)
	\qquad
	\forall i, j \in \{1, \ldots, m\}.
\end{equation}
Then the topological entropy of the interconnected system \eqref{eq:inter} is upper bounded by
\begin{equation}\label{eq:inter-ent-upper-eig-reach}
    h(f, t_0, K) \leq n \pos{\spabs(\hat A^\reachlim_\Network)}.
\end{equation}
\end{thm}

Theorem~\ref{thm:inter-ent-upper-eig-reach} extends \cite[Th.~1]{Liberzon2021} from time-invariant to time-varying systems and replaces bounds over the entire state space with asymptotic bounds over the convex hull of the reachable set.

The upper bounds in \eqref{eq:inter-ent-upper-reach} and \eqref{eq:inter-ent-upper-eig-reach} depend on the asymptotic behavior of the interconnection matrix $ A_\Network(t, x) $ over the convex hull of the reachable set $ \xi(t, t_0, K) $.
As in Section~\ref{sec:tv-ent}, they are obtained by applying different operations to $ A_\Network(t, x) $.

Both bounds hold for any choice of local norms. 
The bound in \eqref{eq:inter-ent-upper-reach} additionally depends on a monotone network norm, whereas the bound in \eqref{eq:inter-ent-upper-eig-reach} does not. 
For fixed choices of local and network norms, either bound may be tighter.
For example, the bound in \eqref{eq:inter-ent-upper-reach} may be tighter because taking the matrix measure before the asymptotic bound can preserve relationships among the entries of the interconnection matrix.
Conversely, the bound in \eqref{eq:inter-ent-upper-eig-reach} may be tighter because the spectral abscissa can be strictly smaller than the matrix measure of the same matrix.
These possibilities are illustrated by Examples~\ref{eg:ltv-ent-upper} and~\ref{eg:lti-ent-upper} in Section~\ref{sec:tv-ent}, where the systems can be viewed as interconnections of two scalar subsystems.
Example~\ref{eg:tv-lorenz} in Section~\ref{sec:eg} further illustrates how the choice of subsystem decomposition affects the interconnected entropy bound, even with fixed local norms.
Additional upper bounds on the entropy of \eqref{eq:inter} are given in Section~\ref{sec:alt}.

\begin{rmk}\label{rmk:inter-cas}
For a cascade (i.e., triangular) interconnection, the interconnection matrix $ A_\Network(t, x) $ is triangular, 
and one can select a triangular comparison matrix $ \hat A^\reachlim_\Network $.
Its spectral abscissa is determined by the diagonal entries.
Thus, the upper bound in \eqref{eq:inter-ent-upper-eig-reach} reduces to
\begin{equation}\label{eq:cas-ent-upper-reach}
	h(f,t_0,K)
	\leq n \pos{\max_{1 \leq i \leq m} \hat\mu^\reachlim_i},
\end{equation}
where
\[
	\hat\mu^\reachlim_i 
	:= \limsup_{t \to \infty} 
	\max_{v \in \co(\xi(t, t_0, K))}
	\mu_{\local i}(\jacob{x_i}{f_i(t, v)}).
\]
The bound in \eqref{eq:cas-ent-upper-reach} depends only on the diagonal Jacobian blocks, although the interconnections may still influence these blocks and the reachable sets over which they are evaluated.
\end{rmk}

For a lower bound on the entropy of \eqref{eq:inter}, note that the trace of the full Jacobian is the sum of the traces of its diagonal blocks.
Thus, the lower bound in \eqref{eq:tv-ent-lower} can also be evaluated using subsystem-level trace bounds, regardless of the interconnection structure.

\section{Alternative upper bounds on entropy}\label{sec:alt}
In this section, we present alternative upper bounds on the topological entropy of general and interconnected nonlinear time-varying systems.

We first establish results complementing Theorems~\ref{thm:tv-ent-bnd} and~\ref{thm:inter-ent-upper-eig-reach}.
Proof sketches are provided in Appendix~\ref{apx:tv-inter-ent-upper-ini}, with complete proofs given in
Appendix~\ref{apx:tv-inter-ent-upper-ini-full}.
\begin{prop}\label{prop:tv-ent-upper-ini}
The topological entropy of the nonlinear time-varying system \eqref{eq:tv} is upper bounded by
\begin{equation}\label{eq:tv-ent-upper-ini}
    h(f, t_0, K) \leq n \pos{\hat\mu^\iniinf},
\end{equation}
where
\begin{equation}\label{eq:tv-meas-lim-ini}
	\hat\mu^\iniinf 
	:= \inf_{t_1 \geq t_0} 
	\limsup_{t \to \infty} 
	\max_{v \in \co(\xi(t_1, t_0, K))} 
	\mu(\jacob{x}{f(t, \xi(t, t_1, v))}).
\end{equation}
\end{prop}

The upper bound in \eqref{eq:tv-ent-upper-ini} depends on the asymptotic behavior of the Jacobian $ \jacob{x}{f(t, x)} $ along solutions starting from the convex hull of the reachable set $ \xi(t_1, t_0, K) $, for arbitrary $ t_1 \geq t_0 $. 
Consequently, it may be tighter than the upper bound in \eqref{eq:tv-ent-upper-reach} if the reachable set $ \xi(t_1, t_0, K) $ is convex for some $ t_1 \geq t_0 $. 

\begin{prop}\label{prop:inter-ent-upper-eig-ini}
Consider the interconnection matrix $ A_\Network(t, x) = [a_{ij}(t, x)] \in \R^{m \times m} $ defined in \eqref{eq:inter-sub-meas}. 
Let $ \bar A^\iniinf_\Network = [\bar a^\iniinf_{ij}] \in \R^{m \times m} $ be any matrix satisfying
\begin{equation}\label{eq:inter-sub-meas-sup-ini}
	\inf_{t_1 \geq t_0} 
	\max_{1 \leq i, j \leq m} \left( 
		\sup_{t \geq t_1} 
		a^\ini_{ij}(t, t_1) - \bar a^\iniinf_{ij} 
	\right) 
	\leq 0,
\end{equation}
where
\[
	a^\ini_{ij}(t, t_1) 
	:= \max_{v \in \co(\xi(t_1, t_0, K))} 
	a_{ij}(t, \xi(t, t_1, v)).
\]
Then the topological entropy of the interconnected system \eqref{eq:inter} is upper bounded by
\begin{equation}\label{eq:inter-ent-upper-eig-ini}
    h(f, t_0, K) \leq n \pos{\spabs(\bar A^\iniinf_\Network)}.
\end{equation}
\end{prop}

The upper bound in \eqref{eq:inter-ent-upper-eig-ini} is obtained by bounding the interconnection matrix along solutions starting from the convex hull of the reachable set at an intermediate time $ t_1 \geq t_0 $. 
Unlike \eqref{eq:inter-sub-meas-lim-reach}, the condition \eqref{eq:inter-sub-meas-sup-ini} involves a supremum over $ t \ge t_1 $, rather than an upper limit as $ t \to \infty $. 
This difference arises from the alternative trajectory-separation bound in Lemma~\ref{lem:inter-soln-upper-ini} in Appendix~\ref{apx:tv-inter-ent-upper-ini}.

We next derive corresponding results complementing Proposition~\ref{prop:tv-ent-upper-eig-reach} and
Corollary~\ref{cor:inter-ent-upper-reach}.
\begin{cor}\label{cor:tv-ent-upper-eig-ini}
Let $ \bar A^\iniinf = [\bar a^\iniinf_{ij}] \in \R^{n \times n} $ be any matrix satisfying
\begin{equation}\label{eq:tv-comp-lim-ini}
	\inf_{t_1 \geq t_0} 
	\max_{1 \leq i, j \leq n} \left( 
		\sup_{t \geq t_1} 
		a^\ini_{ij}(t, t_1) - \bar a^\iniinf_{ij} 
	\right) 
	\leq 0,
\end{equation}
where
\begin{align*}
	a^\ini_{ii}(t, t_1) 
	&:= \max_{v \in \co(\xi(t_1, t_0, K))} 
	\pdv{x_i}{f_i(t, \xi(t, t_1, v))}, \\
	a^\ini_{ij}(t, t_1) 
	&:= \max_{v \in \co(\xi(t_1, t_0, K))} 
	|\pdv{x_j}{f_i(t, \xi(t, t_1, v))}|	
\end{align*}
for diagonal and off-diagonal entries, respectively.
Then the topological entropy of the nonlinear time-varying system \eqref{eq:tv} is upper bounded by
\begin{equation}\label{eq:tv-ent-upper-eig-ini}
    h(f, t_0, K) \leq n \pos{\spabs(\bar A^\iniinf)}.
\end{equation}
\end{cor}

Corollary~\ref{cor:tv-ent-upper-eig-ini} follows from Proposition~\ref{prop:inter-ent-upper-eig-ini} by viewing \eqref{eq:tv} as $ n $ interconnected scalar subsystems, as in the proof of Proposition~\ref{prop:tv-ent-upper-eig-reach}.

\begin{cor}\label{cor:inter-ent-upper-ini}
Consider the interconnection matrix $ A_\Network(t, x) $ defined in \eqref{eq:inter-sub-meas}. 
The topological entropy of the interconnected system \eqref{eq:inter} is upper bounded by
\begin{equation}\label{eq:inter-ent-upper-ini}
    h(f, t_0, K) \leq n \pos{\hat\mu^\iniinf_\Network},
\end{equation}
where
\begin{equation}\label{eq:inter-meas-lim-ini}
	\hat\mu^\iniinf_\Network 
	:= \inf_{t_1 \geq t_0} 
	\limsup_{t \to \infty} 
	\max_{v \in \co(\xi(t_1, t_0, K))} 
	\mu_\Network(A_\Network(t, \xi(t, t_1, v))).
\end{equation}
\end{cor}

Corollary~\ref{cor:inter-ent-upper-ini} follows from Proposition~\ref{prop:tv-ent-upper-ini} and Lemma~\ref{lem:inter-meas-bnd}, as in the proof of Corollary~\ref{cor:inter-ent-upper-reach}.

\section{Numerical illustration}\label{sec:eg}
\begin{eg}\label{eg:tv-lorenz}
We consider a time-varying extension of the classical Lorenz system:
\begin{equation}\label{eq:tv-lorenz}
\begin{aligned}
	\dot x_1 &= \sigma (x_2 - x_1), \\
	\dot x_2 &= \theta(t)\, x_1 - x_2 - x_1 x_3, \\
	\dot x_3 &= x_1 x_2 - \beta x_3,
\end{aligned}
\end{equation}
where $ \theta(t) = \theta_0 + \sin t $ models periodic variation in thermal forcing.
We select $ \sigma = 10 $, $ \beta = 8/3 $, and $ \theta_0 = 28 $, for which the corresponding autonomous Lorenz system exhibits chaotic behavior \cite{Lorenz1963}.
We compare the upper bounds on the entropy of \eqref{eq:tv-lorenz} obtained from Theorems~\ref{thm:tv-ent-bnd} and~\ref{thm:inter-ent-upper-eig-reach}, using different induced norms and subsystem decompositions.

Following a standard quadratic Lyapunov-function construction for the autonomous Lorenz system (see, e.g., \cite{kuznetsov2020lorenz}), define
\[
	V(t, x) = x_1^2 + x_2^2 + (x_3 - \sigma - \theta(t))^2.
\]
Let $ z := x_3 - \sigma - \theta(t) $. 
Differentiating along solutions yields
\begin{align*}
	\dot V
	&= -2 \sigma x_1^2 - 2 x_2^2 - 2 \beta z^2
	- 2 \bigl( \beta (\sigma + \theta(t)) + \dot\theta(t) \bigr) z.
\end{align*}
By Young’s inequality,
\begin{align*}
	\dot V
	&\leq -c V + \frac{\bigl( \beta (\sigma + \theta(t)) + \dot\theta(t) \bigr)^2}{\beta},
\end{align*}
where $ c := \min\{2 \sigma,\, 2,\, \beta\} = 2 $. 
Since
\[
    \sup_{t \geq 0}
    {\bigl| \beta (\sigma + \theta(t)) + \dot\theta(t) \bigr|}
    = \beta (\sigma + \theta_0) + \sqrt{\beta^2 + 1},
\]
the family of time-varying balls
\[
	S_t := \left\{
		x \in \R^3: 
		x_1^2 + x_2^2 + (x_3 - \sigma - \theta(t))^2 
		\leq R^2
	\right\},
\]
where
\[
	R := \frac{
		\beta (\sigma + \theta_0) + \sqrt{\beta^2 + 1}
	}{\sqrt{c \beta}}
\]
is forward invariant.
For any initial set $ K \subset S_0 $, the convexity of $ S_t $ implies
\[
	\co(\xi(t, 0, K)) \subset S_t
	\qquad
	\forall t \geq 0.
\]
Thus, the entropy bounds can be evaluated using the enclosures $ S_t $, as discussed in Remark~\ref{rmk:superset}.

We first evaluate the general upper bound in \eqref{eq:tv-ent-upper-reach} using the induced $ \infty $- and $ 2 $-norms.
Since $ \theta(t) - x_3 = -\sigma - z $, the Jacobian of \eqref{eq:tv-lorenz} is
\[
	\jacob{x}{f(t, x)}
	= \begin{bmatrix}
		-\sigma & \sigma & 0 \\
		-\sigma - z & -1 & -x_1 \\
		x_2 & x_1 & -\beta
	\end{bmatrix}.
\]

For the induced $ \infty $-norm, the corresponding matrix measure satisfies
\begin{align*}
	&\mu_\infty(\jacob{x}{f(t, x)}) \\
	&\qquad = \max \left\{
		0,\,
		-1 + |\sigma + z| + |x_1|,\,
		-\beta + |x_2| + |x_1|
	\right\}.
\end{align*}
Since $ x_1^2 + x_2^2 + z^2 \leq R^2 $ on $ S_t $ and $ \sigma + \beta > 1 $,
\[
	\sup_{t \geq 0} 
	\max_{x \in S_t} 
	\mu_\infty(\jacob{x}{f(t, x)})
	= \sigma - 1 + \sqrt{2} R.
\]
Thus, Theorem~\ref{thm:tv-ent-bnd} yields
\[
	h(f, 0, K) \leq 3 (\sigma - 1 + \sqrt{2} R).
\]

For the induced $ 2 $-norm, the corresponding matrix measure $ \mu_2(\jacob{x}{f(t, x)}) $ is the largest eigenvalue of
\[
	\frac{\jacob{x}{f(t, x)} + \jacob{x}{f(t, x)}^\top}{2}
	= \begin{bmatrix}
		-\sigma & -z/2 & x_2/2 \\
		-z/2 & -1 & 0 \\
		x_2/2 & 0 & -\beta
\end{bmatrix},
\]
where the skew-symmetric coupling between $ x_2 $ and $ x_3 $ does not contribute to the symmetric part.
Since $ \beta > 1 $, the largest eigenvalue over $ S_t $ is maximized at $ |z| = R $ and $ x_2 = 0 $, giving
\[
	\sup_{t \geq 0} 
	\max_{x \in S_t} 
	\mu_2(\jacob{x}{f(t, x)})
	= \frac{-(\sigma + 1) + \sqrt{(\sigma - 1)^2 + R^2}}{2}.
\]
Thus, Theorem~\ref{thm:tv-ent-bnd} yields
\[
	h(f, 0, K)
	\leq \frac{3}{2} \left(
		-(\sigma+1)
		+\sqrt{(\sigma-1)^2+R^2}
	\right).
\]

We next evaluate the interconnected upper bound in \eqref{eq:inter-ent-upper-eig-reach} using three subsystem decompositions.
For each decomposition, we construct a Metzler comparison matrix satisfying \eqref{eq:inter-sub-meas-lim-reach} with respect to the interconnection matrix defined by \eqref{eq:inter-sub-meas}.

First, viewing \eqref{eq:tv-lorenz} as three interconnected scalar subsystems and using $ |x_1|, |x_2|, |z| \leq R $ on $ S_t $, we obtain the comparison matrix
\[
	\hat A_{\Network}^{(1)}
	= \begin{bmatrix}
		-\sigma & \sigma & 0 \\
		\sigma + R & -1 & R \\
		R & R & -\beta
	\end{bmatrix}.
\]
Thus, Theorem~\ref{thm:inter-ent-upper-eig-reach} yields
\[
	h(f, 0, K)
	\leq 3 \pos{\spabs(\hat A_{\Network}^{(1)})}.
\]

Second, consider two interconnected subsystems with states $ y_1 = (x_1, x_2) $ and $ y_2 = x_3 $, using induced $ 2 $-norms locally.
Since the symmetric part of the first diagonal Jacobian block is
\[
    \begin{bmatrix}
        -\sigma & -z/2 \\
        -z/2 & -1
    \end{bmatrix},
\]
the interconnection matrix satisfies
\[
    a_{11}(t, x) 
    \leq \frac{-(\sigma + 1) + \sqrt{(\sigma - 1)^2 + R^2}}{2}
    =: \psi,
\]
and
\[
    a_{12}(t, x) \leq R,
    \qquad
    a_{21}(t, x) \leq R,
    \qquad
    a_{22}(t, x) = -\beta.
\]
Hence the corresponding comparison matrix is
\[
	\hat A_\Network^{(2)}
	= \begin{bmatrix}
		\psi & R \\
		R & -\beta
	\end{bmatrix}.
\]
Thus, Theorem~\ref{thm:inter-ent-upper-eig-reach} yields
\begin{align*}
    h(f, 0, K)
    &\leq 
    \frac{3}{2}
    \left(
        \psi - \beta + \sqrt{(\psi + \beta)^2 + 4 R^2}
    \right).
\end{align*}

Finally, consider two interconnected subsystems with states $ y_1 = x_1 $ and $ y_2 = (x_2, x_3) $, again using induced $ 2 $-norms locally.
Since the symmetric part of the second diagonal Jacobian block is
\[
	\begin{bmatrix} 
		-1 & 0 \\ 
		0 & -\beta 
	\end{bmatrix},
\]
the interconnection matrix satisfies
\begin{align*}
	a_{11}(t, x) &= -\sigma, 
	&a_{12}(t, x) &= \sigma, \\
	a_{21}(t, x) &\leq \sigma + R,
	&a_{22}(t, x) &= -1.
\end{align*}
Hence the corresponding comparison matrix is
\[
	\hat A_{\Network}^{(3)}
	= \begin{bmatrix}
		-\sigma & \sigma \\
		\sigma + R & -1
	\end{bmatrix}.
\]
Thus, Theorem~\ref{thm:inter-ent-upper-eig-reach} yields
\begin{align*}
	h(f, 0, K)
	&\leq 
	\frac{3}{2}
    \left(
        -(\sigma + 1) + \sqrt{(\sigma - 1)^2 + 4 \sigma (\sigma + R)}
    \right).
\end{align*}

The resulting numerical bounds are summarized in Table~\ref{tab:eg}.
For the selected parameters, $ R \approx 45.11 $.
\begin{table}[!ht]
\renewcommand{\arraystretch}{1.2}
\caption{Upper bounds on the entropy of the time-varying Lorenz system.}
\label{tab:eg}
\centering
\begin{tabular}{cc}
\hline
Method & Upper bound \\
\hline
\eqref{eq:tv-ent-upper-reach}, induced $ \infty $-norm
& $ 218.39 $ \\
\eqref{eq:tv-ent-upper-reach}, induced $ 2 $-norm 	
& $ 52.50 $ \\
\eqref{eq:inter-ent-upper-eig-reach}, scalar subsystems 	
& $ 154.71 $ \\
\eqref{eq:inter-ent-upper-eig-reach}, $ (x_1, x_2) $--$ x_3 $, induced $ 2 $-norm 
& $ 160.93 $ \\
\eqref{eq:inter-ent-upper-eig-reach}, $ x_1 $--$ (x_2, x_3) $, induced $ 2 $-norm 	
& $ 55.21 $ \\
\hline
\end{tabular}\end{table}

The comparison illustrates the effects of norm selection and subsystem decomposition on the entropy bounds.
The general bound \eqref{eq:tv-ent-upper-reach} obtained using the induced $ 2 $-norm is substantially tighter than that obtained using the induced $ \infty $-norm.
The interconnected bound \eqref{eq:inter-ent-upper-eig-reach} with the $ x_1 $--$ (x_2, x_3) $ decomposition also improves substantially upon the general induced $ \infty $-norm bound, although it is slightly larger than the general induced $ 2 $-norm bound.

The subsystem decomposition significantly affects the interconnected bound \eqref{eq:inter-ent-upper-eig-reach}, even with fixed local norms.
In particular, grouping $ x_2 $ and $ x_3 $ into one subsystem preserves their skew-symmetric coupling within a diagonal Jacobian block, where it does not contribute to the matrix measure.
In contrast, separating these variables introduces nonnegative off-diagonal block norms and loses this cancellation.
Consequently, the $ x_1 $--$ (x_2, x_3) $ decomposition yields a substantially tighter bound than the other decompositions considered.

Finally, the example illustrates the usefulness of the time-varying enclosures discussed in Remark~\ref{rmk:superset}.
Centering $ S_t $ at $ x_3 = \sigma + \theta(t) $ simplifies the Jacobian entry $ \theta(t) - x_3 $ to $ -\sigma - z $, allowing the time dependence to be incorporated directly into the entropy bounds.
Also, using the convex enclosures $ S_t $, Propositions~\ref{prop:tv-ent-upper-ini} and~\ref{prop:inter-ent-upper-eig-ini} yield the same numerical bounds as Theorems~\ref{thm:tv-ent-bnd} and~\ref{thm:inter-ent-upper-eig-reach}, respectively.
\end{eg}

\section{Proofs of main results}\label{sec:proof}
In this section, we prove the main results.
We first establish a bound on trajectory separation and construct universal spanning sets.
We then combine these results to prove Theorem~\ref{thm:inter-ent-upper-eig-reach}, followed by Theorem~\ref{thm:tv-ent-bnd} and Proposition~\ref{prop:tv-ent-upper-eig-reach}.

\subsection{Subsystem-wise trajectory-separation bound}\label{ssec:soln-vol}
We establish a subsystem-wise upper bound on the separation between trajectories of the interconnected system \eqref{eq:inter}. 
For each $ i \in \{1, \ldots, m\} $, let $ \xi_i(t, t_0, x) \in \R^{n_i} $ denote the state of the $ i $-th subsystem at time $ t $ along the solution of \eqref{eq:inter} with initial state $ x $ at time $ t_0 $. 
The overall solution is then given by $ \xi(t, t_0, x) = (\xi_1(t, t_0, x), \ldots, \xi_m(t, t_0, x)) $.

\begin{lem}\label{lem:inter-soln-upper-reach}
Consider the interconnection matrix $ A_\Network(t, x) = [a_{ij}(t, x)] \in \R^{m \times m} $ defined in \eqref{eq:inter-sub-meas}.
Let $ \bar A_\Network^\reach(t, t_0, K) = [\bar a^\reach_{ij}(t, t_0, K)] \in \R^{m \times m} $ be any matrix-valued function satisfying
\begin{equation}\label{eq:inter-sub-meas-max-reach}
	\bar a^\reach_{ij}(t, t_0, K) 
	\geq \sup_{s \in [t_0, t]} 
	\max_{v \in \co(\xi(s, t_0, K))} a_{ij}(s, v)
\end{equation}
for all $ t \geq t_0 $ and $ i, j \in \{1, \ldots, m\} $.
Then for all initial states $ x = (x_1, \ldots, x_m), \bar x = (\bar x_1, \ldots, \bar x_m) \in K $ with $ x_i, \bar x_i \in \R^{n_i} $, the solutions of the interconnected system \eqref{eq:inter} satisfy
\begin{align}
	&\begin{bmatrix}
		|\xi_1(t, t_0, \bar x) - \xi_1(t, t_0, x)|_{\local 1} \\
		\vdots \\
		|\xi_m(t, t_0, \bar x) - \xi_m(t, t_0, x)|_{\local m}
	\end{bmatrix} \notag\\
	&\qquad \leq e^{\bar A_\Network^\reach(t, t_0, K) (t - t_0)} 
	\begin{bmatrix}
		|\bar x_1 - x_1|_{\local 1} \\
		\vdots \\
		|\bar x_m - x_m|_{\local m}
	\end{bmatrix}
	\qquad
	\forall t \geq t_0. \label{eq:inter-soln-upper-reach}
\end{align}
\end{lem}

Lemma~\ref{lem:inter-soln-upper-reach} extends \cite[Prop.~1]{ArcakMaidens2018} by introducing the matrix-valued function $ \bar A_\Network^\reach(t, t_0, K) $ in place of a constant matrix. 
This extension is used to derive the upper bound on entropy in Theorem~\ref{thm:inter-ent-upper-eig-reach}. 
The proof, given in Appendix~\ref{apx:inter-soln-upper-reach}, is inspired by the variational construction in \cite[Prop.~1]{ArcakMaidens2018}, which builds on earlier results including \cite[Th.~1]{Sontag2010} and \cite[Th.~4.2]{BoichenkoLeonov1998}.

\subsection{Universal spanning sets}\label{ssec:ent-span}
We construct spanning sets by extending a grid-based approach from \cite{YangSchmidtLiberzonHespanha2020}. 
This construction is independent of the system dynamics.

As $ K $ is compact, there exists a closed hypercube $ B $ of radius $ r > 0 $ centered at the origin such that $ K \subset B $.
For a vector $ \theta = (\theta_1, \ldots, \theta_n) \in \R_{>0}^n $, which may depend on $ T $ and $ \varepsilon $, define the following grid in $ B $:
\begin{equation}\label{eq:grid-dfn}
    G(\theta) 
    := \{(k_1 \theta_1, \ldots, k_n \theta_n) \in B: k_1, \ldots, k_n \in \Z\}.
\end{equation}
Let $ R(x, \theta) $ denote the open hyperrectangle in $ \R^n $ centered at $ x = (x_1, \ldots, x_n) $ with side lengths $ 2 \theta_1, \ldots, 2 \theta_n $, that is,
\begin{align}
    R(x, \theta) &:= \{(\bar x_1, \ldots, \bar x_n) \in \R^n: \notag\\
    &\qquad |\bar x_1 - x_1| < \theta_1, \ldots, |\bar x_n - x_n| < \theta_n\}. \label{eq:grid-rect}
\end{align}
Then
\[
    B \subset \bigcup_{x \in G(\theta)} R(x, \theta).
\]
Next, for each $ x \in G(\theta) $ such that $ R(x, \theta) \cap K \neq \emptyset $, select one point $ \bar x \in R(x, \theta) \cap K $, and let $ E(\theta) $ be the set of the selected points.
Then
\[
	K \subset \bigcup_{x \in E(\theta)} R(x, 2 \theta)
\]
by the triangle inequality, and the cardinality of $ E(\theta) $ is upper bounded by
\begin{equation*}
	\card{E(\theta)} 
	\leq \card{G(\theta)} 
	\leq \prod_{i=1}^{n} \left( 
		2 \left\lfloor \frac{r}{\theta_i} \right\rfloor + 1 
	\right).
\end{equation*}

Given a time horizon $ T \geq 0 $ and a radius $ \varepsilon > 0 $, by comparing the hyperrectangle $ R(x, 2 \theta) $ to the open ball $ B_{f, t_0}(x, \varepsilon, T) $ defined in \eqref{eq:ball-dfn}, we obtain the following lemma.
The proof is given in Appendix~\ref{apx:grid-span}.
\begin{lem}\label{lem:grid-span}
If the vector $ \theta $ is selected so that
\begin{equation}\label{eq:grid-span-cond}
	(R(x, 2 \theta) \cap K) 
	\subset B_{f, t_0}(x, \varepsilon, T) 
	\qquad 
	\forall x \in E(\theta),
\end{equation}
then the set $ E(\theta) $ is $ (T, \varepsilon) $-spanning. 
Moreover, if \eqref{eq:grid-span-cond} holds for all $ T > 0 $ and $ \varepsilon > 0 $, and
\begin{equation}\label{eq:grid-span-ent-cond}
    \limsup_{\varepsilon \searrow 0} 
    \limsup_{T \to \infty} 
    \frac{\log\theta_i}{T} 
    \leq 0 
    \qquad 
    \forall i \in \{1, \ldots, n\},
\end{equation}
then the topological entropy of \eqref{eq:tv} is upper bounded by
\begin{equation}\label{eq:grid-span-ent}
    h(f, t_0, K) 
    \leq \limsup_{\varepsilon \searrow 0} 
    \limsup_{T \to \infty} 
    \sum_{i=1}^{n} \frac{\log(1/\theta_i)}{T}.
\end{equation}
\end{lem}

Note that \eqref{eq:grid-span-ent-cond} holds if all $ \theta_i $ are nonincreasing in $ T $.

Lemma~\ref{lem:grid-span} extends the results on spanning sets and upper bounds on entropy from \cite[Lemma~2]{YangSchmidtLiberzonHespanha2020} and \cite[Lemma~6]{YangLiberzonHespanha2023}, from switched linear and nonlinear systems, respectively, to general nonlinear time-varying systems.

\subsection{Proof of Theorem~\ref{thm:inter-ent-upper-eig-reach}}\label{ssec:inter-ent-upper-eig-reach-pf}
Combining the trajectory-separation bound in Lemma~\ref{lem:inter-soln-upper-reach} with the spanning set construction in Lemma~\ref{lem:grid-span} yields the following intermediate upper bound on entropy.
The proof is given in Appendix~\ref{apx:inter-ent-upper-eig-reach}.
\begin{lem}\label{lem:inter-ent-upper-eig-reach}
Consider the interconnection matrix $ A_\Network(t, x) = [a_{ij}(t, x)] \in \R^{m \times m} $ defined in \eqref{eq:inter-sub-meas}.
Let $ \bar A_\Network^\reach(t, t_0, K) = [\bar a^\reach_{ij}(t, t_0, K)] \in \R^{m \times m} $ be any locally bounded matrix-valued function satisfying \eqref{eq:inter-sub-meas-max-reach} for all $ t \geq t_0 $ and $ i, j \in \{1, \ldots, m\} $.
Then the topological entropy of \eqref{eq:inter} is upper bounded by
\begin{equation}\label{eq:inter-ent-upper-eig-reach-temp}
    h(f, t_0, K) \leq \limsup_{T \to \infty} \frac{n}{T} \sup_{t \in [t_0, t_0 + T]} \log\|e^{\bar A_\Network^\reach(t, t_0, K) (t - t_0)}\|_\Network.
\end{equation}
\end{lem}

We now establish the upper bound on entropy in Theorem~\ref{thm:inter-ent-upper-eig-reach}.
\begin{proof}[Proof of Theorem~\ref{thm:inter-ent-upper-eig-reach}]
By the upper limit in \eqref{eq:inter-sub-meas-lim-reach}, for each $ \delta_1 > 0 $, there exists a sufficiently large $ t_1 \geq t_0 $ such that 
\begin{equation*}
	\sup_{t \geq t_1} 
	\max_{v \in \co(\xi(t, t_0, K))} 
	a_{ij}(t, v) 
	\leq \hat a^\reachlim_{ij} + \delta_1 
	\qquad 
	\forall i, j \in \{1, \ldots, m\}.
\end{equation*}

Consider the interconnected system \eqref{eq:inter} with the new initial time $ t_1 $ and the corresponding initial set $ K_1 := \xi(t_1, t_0, K) $.
As $ \xi(t, t_0, K) = \xi(t, t_1, K_1) $, the matrix-valued function $ \bar A_\Network^\reach(t, t_1, K_1) := \hat A^\reachlim_\Network + \delta_1 \oneMx $, where $ \oneMx $ denotes the matrix of ones in $ \R^{m \times m} $, satisfies \eqref{eq:inter-sub-meas-max-reach} for all $ t \geq t_1 $.
Thus, applying Lemma~\ref{lem:inter-ent-upper-eig-reach} with the new initial time $ t_1 $ and initial set $ K_1 $ yields
\[
    h(f, t_1, K_1) \leq \limsup_{T \to \infty} \frac{n}{T} \max_{t \in [t_1, t_1 + T]} \log\|e^{(\hat A^\reachlim_\Network + \delta_1 \oneMx) (t - t_1)}\|_\Network.
\]
Therefore, by Lemma~\ref{lem:ini-time},
\begin{align*}
    h(f, t_0, K)
    &= h(f, t_1, K_1) \\
    &\leq \limsup_{T \to \infty} \frac{n}{T} 
    \max_{t \in [0, T]} 
    \log\|e^{(\hat A^\reachlim_\Network + \delta_1 \oneMx) t}\|_\Network.
\end{align*}

By characterizing the spectral abscissa through the asymptotic growth rate of the matrix exponential \cite[Fact~11.15.7, p.~690]{Bernstein2009}, we obtain
\[
	\spabs(\hat A^\reachlim_\Network + \delta_1 \oneMx) 
	= \lim_{t \to \infty} 
	\frac{1}{t} \log\|e^{(\hat A^\reachlim_\Network + \delta_1 \oneMx) t}\|_\Network.
\]
Thus, for each $ \delta_2 > 0 $, there exists a sufficiently large $ t_2 \geq 0 $ such that
\[
	\log\|e^{(\hat A^\reachlim_\Network + \delta_1 \oneMx) t}\|_\Network 
	\leq (\spabs(\hat A^\reachlim_\Network + \delta_1 \oneMx) + \delta_2)\, t \qquad \forall t \geq t_2.
\]
Consequently, for all $ T > t_2 $, we have
\[
	\max_{t \in [0, T]} 
	\log\|e^{(\hat A^\reachlim_\Network + \delta_1 \oneMx) t}\|_\Network 
	\leq \max \left\{ 
		(\spabs(\hat A^\reachlim_\Network + \delta_1 \oneMx) + \delta_2)\, T,\, 
		C_2 
	\right\},
\]
where $ C_2 := \max_{t \in [0, t_2]} \log\|e^{(\hat A^\reachlim_\Network + \delta_1 \oneMx) t}\|_\Network $ is independent of $ T $ and thus satisfies $ C_2/T \to 0 $ as $ T \to \infty $. Therefore,
\begin{align*}
    h(f, t_0, K) 
    &\leq n \max \left\{ \spabs(\hat A^\reachlim_\Network + \delta_1 \oneMx) + \delta_2,\, 0 \right\}.
\end{align*}
Since $ \delta_1, \delta_2 > 0 $ are arbitrary and the spectral abscissa is continuous with respect to the matrix entries (see, e.g., \cite[Sec.~2.4.9]{HornJohnson2013}), the upper bound in \eqref{eq:inter-ent-upper-eig-reach} follows.
\end{proof}

\subsection{Proof of Theorem~\ref{thm:tv-ent-bnd} and Proposition~\ref{prop:tv-ent-upper-eig-reach}}\label{ssec:tv-ent-bnd-pf}
The upper bound in Theorem~\ref{thm:tv-ent-bnd} and Proposition~\ref{prop:tv-ent-upper-eig-reach} follow from Theorem~\ref{thm:inter-ent-upper-eig-reach}.
To establish the lower bound in Theorem~\ref{thm:tv-ent-bnd}, we first derive the following volume bound.
The proof is given in Appendix~\ref{apx:tv-vol-lower}.
\begin{lem}\label{lem:tv-vol-lower}
The reachable set of \eqref{eq:tv} starting from $ K $ satisfies
\begin{equation}\label{eq:tv-vol-lower}
    \vol(\xi(t, t_0, K)) \geq e^{\gamma(t, t_0)} \vol(K) \qquad \forall t \geq t_0,
\end{equation}
where
\begin{equation*}\label{eq:tv-tr-int}
    \gamma(t, t_0) := \min_{v \in K} \int_{t_0}^{t} \tr(\jacob{x}{f(s, \xi(s, t_0, v))}) \dl{s}.
\end{equation*}
\end{lem}

We now establish the upper and lower bounds on entropy in Theorem~\ref{thm:tv-ent-bnd}.
\begin{proof}[Proof of Theorem~\ref{thm:tv-ent-bnd}]
First, the general system \eqref{eq:tv} can be viewed as an interconnected system \eqref{eq:inter} with a single subsystem. 
Thus, the upper bound in \eqref{eq:tv-ent-upper-reach} follows from \eqref{eq:inter-ent-upper-eig-reach} by setting $ m = 1 $ and letting the local matrix measure $ \mu_{\local 1}(\cdot) = \mu(\cdot) $.

We next establish the lower bound in \(\eqref{eq:tv-ent-lower}\) using a volume argument. 
For arbitrary $ T \geq 0 $ and $ \varepsilon > 0 $, the lower bound in \eqref{eq:tv-vol-lower} on the volume of the reachable set implies
\begin{equation*}
    \vol(\xi(t_0 + T, t_0, K)) 
    \geq e^{\gamma(t_0 + T, t_0)} \vol(K).
\end{equation*}
Let $ E $ be a minimal $ (T, \varepsilon) $-spanning set. By the equivalence of norms on finite-dimensional vector spaces, there exists a constant $ r_\infty > 0 $ such that
\begin{equation*}
	|v|_\infty \leq r_\infty |v| 
	\qquad 
	\forall v \in \R^n
\end{equation*}
for the norm $ |\cdot| $. 
Thus,
\begin{align*}
	&\quad\,\, \xi(t_0 + T, t_0, K) \\
	&\subset \bigcup_{x \in E} \{
		\bar x \in \R^n: 
		|\bar x - \xi(t_0 + T, t_0, x)| < \varepsilon
	\} \\
	&\subset \bigcup_{x \in E} \{
		\bar x \in \R^n: 
		|\bar x - \xi(t_0 + T, t_0, x)|_\infty < r_\infty \varepsilon
	\},
\end{align*}
and hence
\begin{align*}
    &\quad\,\, \vol(\xi(t_0 + T, t_0, K)) \\
    &\leq \sum_{x \in E} \vol(\{\bar x \in \R^n: |\bar x - \xi(t_0 + T, t_0, x)|_\infty < r_\infty \varepsilon\}) \\
    &= (2 r_\infty \varepsilon)^n \card{E}.
\end{align*}
Therefore, the minimal cardinality of a $ (T, \varepsilon) $-spanning set satisfies
\begin{align*}
    S(f, t_0, \varepsilon, T, K) 
    = \card{E} 
    &\geq \frac{\vol(\xi(t_0 + T, t_0, K))}{(2 r_\infty \varepsilon)^n} \\
    &\geq \frac{e^{\gamma(t_0 + T, t_0)} \vol(K)}{(2 r_\infty \varepsilon)^n}.
\end{align*}
Consequently, by the definition of entropy \eqref{eq:ent-dfn}, we have
\begin{align*}
    &\quad\,\, h(f, t_0, K) \\
    &\geq \liminf_{\varepsilon \searrow 0} 
    \limsup_{T \to \infty} \frac{
    	\gamma(t_0 + T, t_0) 
    	+ \log(\vol(K)/(2 r_\infty \varepsilon)^n)
    }{T}  \\
    &= \limsup_{T \to \infty} \frac{\gamma(t_0 + T, t_0)}{T},
\end{align*}
where we use the fact that $ C_\varepsilon := \log(\vol(K)/(2 r_\infty \varepsilon)^n) $ is independent of $ T $ and thus satisfies $ C_\varepsilon/T \to 0 $ as $ T \to \infty $.
Therefore,
\begin{align*}
    &\quad\,\, h(f, t_0, K) \\
    &\geq \limsup_{T \to \infty} \frac{1}{T} \min_{v \in K} 
    \int_{t_0}^{t_0 + T} \tr(\jacob{x}{f(s, \xi(s, t_0, v))}) \dl{s} \\
    &\geq \limsup_{T \to \infty} \frac{1}{T} 
    \int_{t_0}^{t_0 + T} \left( \min_{v \in K} \tr(\jacob{x}{f(s, \xi(s, t_0, v))}) \right) \dl{s}.
\end{align*}

By the lower limit in \eqref{eq:tv-tr-lim}, for each $ \delta > 0 $, there exists a sufficiently large $ t_\delta \geq t_0 $ such that 
\begin{equation*}
	\min_{v \in K} \tr(\jacob{x}{f(t, \xi(t, t_0, v))}) 
	\geq \check\chi - \delta
	\qquad 
	\forall t \geq t_\delta.
\end{equation*}
Consequently,
\[
	\int_{t_0}^{t_0 + T} \left( 
		\min_{v \in K} \tr(\jacob{x}{f(s, \xi(s, t_0, v))}) 
	\right) \dl{s} 
	\geq (\check\chi - \delta) T + C_\delta
\]
for all $ T \geq t_\delta - t_0 $, where 
\[
	C_\delta 
	:= \int_{t_0}^{t_\delta} \left( 
		\min_{v \in K} \tr(\jacob{x}{f(s, \xi(s, t_0, v))}) 
		- (\check\chi - \delta) 
	\right) \dl{s}
\]
is independent of $ T $ and thus satisfies $ C_\delta/T \to 0 $ as $ T \to \infty $. 
Therefore,
\[
    h(f, t_0, K) \geq \check\chi - \delta.
\]
Since $ \delta > 0 $ is arbitrary and $ h(f, t_0, K) \geq 0 $, the lower bound in \eqref{eq:tv-ent-lower} follows.
\end{proof}

\begin{proof}[Proof of Proposition~\ref{prop:tv-ent-upper-eig-reach}]
The general system \eqref{eq:tv} can be viewed as an interconnected system \eqref{eq:inter} composed of $ n $ scalar subsystems. 
Thus, the upper bound in \eqref{eq:tv-ent-upper-eig-reach} follows from \eqref{eq:inter-ent-upper-eig-reach} by setting $ m = n $ and letting the local norms be the absolute value (in which case the local matrix measure of a scalar is the scalar itself).
\end{proof}

\section{Conclusion}\label{sec:end}
This paper studied the topological entropy of general and interconnected nonlinear time-varying systems.
We examined the dependence of entropy on the initial time and initial set, and established general upper and lower bounds using asymptotic properties of the system Jacobian.
For interconnected systems, we derived upper bounds based on the matrix measure of an interconnection
matrix and the spectral abscissa of a Metzler comparison matrix.
The latter incorporates interconnection structure through subsystem-wise trajectory-separation analysis. 
It recovers the two general upper bounds by viewing the system as an interconnection of a single subsystem or scalar subsystems, respectively.

The proposed approach has several notable features:
\begin{itemize}
	\item It allows enclosures of reachable sets to be used to obtain more tractable bounds, potentially at the expense of tightness.

	\item It separates trajectory-separation analysis, spanning-set construction, and asymptotic entropy estimates, providing flexibility in adapting these components to different classes of systems.

	\item It supports arbitrary subsystem decompositions and facilitates modular analysis, allowing bounds on the Jacobian blocks to be obtained separately and combined to estimate the entropy of the overall system.
		The choice of norm and decomposition can potentially be optimized based on the system structure.
\end{itemize}

Future work will investigate improvements that further exploit time-varying dynamics, extensions to special interconnection structures such as cascades, and the computational complexity of the proposed bounds.

\appendix
\section{Proofs of technical lemmas}
\subsection{Proof of Lemma~\ref{lem:ini-time}}\label{apx:ini-time}
\pushQED{\qed}
Fix $ t_1 \geq t_0 $ and let $ K_1 := \xi(t_1, t_0, K) $. 
Consider arbitrary $ T > t_1 - t_0 $ and $ \varepsilon > 0 $.

First, let $ E_0 $ be a minimal $ (T, \varepsilon) $-spanning set for $ K $ at $ t_0 $, and $ E_1 := \xi(t_1, t_0, E_0) $. 
We show that $ E_1 $ is a $ (t_0 + T - t_1, \varepsilon) $-spanning set for $ K_1 $ at $ t_1 $ using the composition property of solutions.
For any $ z \in K_1 $, there exists $ x \in K $ such that $ z = \xi(t_1, t_0, x) $.
By the spanning properties of $ E_0 $, there exists $ \bar x \in E_0 $ such that $ x \in B_{f, t_0}(\bar x, \varepsilon, T) $. 
Then $ \bar z := \xi(t_1, t_0, \bar x) \in E_1 $ satisfies
\[
	|\xi(t, t_1, z) - \xi(t, t_1, \bar z)| 
	= |\xi(t, t_0, x) - \xi(t, t_0, \bar x)| 
	< \varepsilon
\]
for all $ t \in [t_1, t_0 + T] $, that is, $ z \in B_{f, t_1}(\bar z, \varepsilon, t_0 + T - t_1) $. 
Thus,
\[
	K_1 \subset \cup_{\bar z \in E_1} B_{f, t_1}(\bar z, \varepsilon, t_0 + T - t_1).
\] 
Therefore, 
\begin{equation*}
	S(f, t_1, \varepsilon, t_0 + T - t_1, K_1) 
	\leq \card{E_0} 
	= S(f, t_0, \varepsilon, T, K).
\end{equation*}
Since $ T/(t_0 + T - t_1) \to 1 $ as $ T \to \infty $, the definition of entropy \eqref{eq:ent-dfn} yields $ h(f, t_1, K_1) \leq h(f, t_0, K) $.%

Conversely, let $ \bar E_0 $ be a minimal $ (t_1 - t_0, \varepsilon) $-spanning set for $ K $ at $ t_0 $, and $ \bar E_1 $ a minimal $ (t_0 + T - t_1, \varepsilon) $-spanning set for $ K_1 $ at $ t_1 $.
For each pair $ (\bar x, \bar z) \in \bar E_0 \times \bar E_1 $, let
\begin{align*}
	C_{\bar x, \bar z}
	&:= K \cap B_{f, t_0}(\bar x, \varepsilon, t_1 - t_0) \\
	&\;\quad\, \cap \{x: \xi(t_1, t_0, x) \in B_{f, t_1}(\bar z, \varepsilon, t_0 + T - t_1)\}.
\end{align*}
For every nonempty $ C_{\bar x, \bar z} $, select one point $ \tilde x \in C_{\bar x, \bar z} $, and let $ E $ be the set of the selected points. Then $ \card{E} \leq \card{\bar E_0}\, \card{\bar E_1} $.
We show that $ E $ is a $ (T, 2 \varepsilon) $-spanning set for $ K $ at $ t_0 $.
By the spanning properties of $ \bar E_0 $ and $ \bar E_1 $, for any $ x \in K $, there exist $ \bar x \in \bar E_0 $ and $ \bar z \in \bar E_1 $ such that
\[
	x \in B_{f, t_0}(\bar x, \varepsilon, t_1 - t_0),
	\qquad
	\xi(t_1, t_0, x) \in B_{f, t_1}(\bar z, \varepsilon, t_0 + T - t_1).
\]
Hence $ x \in C_{\bar x, \bar z} $, and by construction, there exists $ \tilde x \in E \cap C_{\bar x, \bar z} $.
Applying the triangle inequality with the intermediate points $ \xi(t, t_0, \bar x) $ for $ t \in [t_0, t_1] $ and $ \xi(t, t_1, \bar z) $ for $ t \in [t_1, t_0 + T] $, and using the composition property of solutions, we obtain%
\[
	|\xi(t, t_0, x) - \xi(t, t_0, \tilde x)| 
	< 2 \varepsilon
	\qquad 
	\forall t \in [t_0, t_0 + T].
\]
Therefore,
\begin{align*}
	S(f, t_0, 2\varepsilon, T, K) 
	&\leq S(f, t_0, \varepsilon, t_1 - t_0, K) \\
	&\quad \,\times S(f, t_1, \varepsilon, t_0 + T - t_1, K_1).
\end{align*}
Since the first factor is independent of $ T $, the definition of entropy \eqref{eq:ent-dfn} yields $ h(f, t_0, K) \leq h(f, t_1, K_1) $.

Combining the two inequalities completes the proof.
\popQED

\subsection{Proof of Lemma~\ref{lem:inter-soln-upper-reach}}\label{apx:inter-soln-upper-reach}
We first establish a subsystem-wise trajectory bound for the \emph{interconnected LTV} system 
\begin{equation}\label{eq:inter-ltv}
	\dot x_i = \sum_{j=1}^{m} A_{ij}(t)\, x_j, \quad t \geq t_0, \qquad i \in \{1, \ldots, m\},
\end{equation}
where $ A_{ij}: \R \to \R^{n_i \times n_j} $ are piecewise-continuous matrix-valued functions.

\begin{lem}\label{lem:inter-ltv-soln-upper}
Let $ \bar A^\lin_\Network(t, t_0) = [\bar a^\lin_{ij}(t, t_0)] \in \R^{m \times m} $ be any matrix-valued function satisfying
\begin{equation}\label{eq:inter-ltv-meas-max}
\begin{aligned}
	\bar a^\lin_{ii}(t, t_0) 
	&\geq \sup_{s \in [t_0, t]} 
	\mu_{\local i}(A_{ii}(s)), \\
	\bar a^\lin_{ij}(t, t_0) 
	&\geq \sup_{s \in [t_0, t]} 
	\|A_{ij}(s)\|_{\local ij}
\end{aligned} 
\qquad
\forall t \geq t_0
\end{equation}
for all diagonal and off-diagonal entries, respectively.
For all initial states $ x = (x_1, \ldots, x_m) \in \R^n $ with $ x_i \in \R^{n_i} $, the solution of the interconnected LTV system \eqref{eq:inter-ltv} satisfies
\begin{equation}\label{eq:inter-ltv-soln-upper}
	\begin{bmatrix}
		|\xi_1(t_1, t_0, x)|_{\local 1} \\
		\vdots \\
		|\xi_m(t_1, t_0, x)|_{\local m}
	\end{bmatrix} 
	\leq e^{\bar A^\lin_\Network(t_1, t_0) (t_1 - t_0)} 
	\begin{bmatrix}
		|x_1|_{\local 1} \\
		\vdots \\
		|x_m|_{\local m}
	\end{bmatrix} \qquad \forall t_1 \geq t_0.
\end{equation}
\end{lem}

\begin{proof}
Fix $ t_1 \geq t_0 $.
For each $ i \in \{1, \ldots, m\} $, the right-sided derivative of the map $ t \mapsto |\xi_i(t, t_0, x)|_{\local i} $, denoted by $ \rpdv{t}{|\xi_i(t, t_0, x)|_{\local i}} $, satisfies the following inequalities at every continuity point of $ A(t) $ in $ [t_0, t_1) $:
\begin{align*}
	&\quad\,\, \rpdv{t}{|\xi_i(t, t_0, x)|_{\local i}} \\
	&= \lim_{s \searrow 0} 
	\frac{
		|\xi_i(t + s, t_0, x)|_{\local i} - |\xi_i(t, t_0, x)|_{\local i}
	}{s} \\
	&= \lim_{s \searrow 0} 
	\frac{
		|\xi_i(t, t_0, x) + s \sum_{j=1}^{m} A_{ij}(t)\, \xi_j(t, t_0, x)|_{\local i} - |\xi_i(t, t_0, x)|_{\local i}
	}{s} \\
	&\leq \lim_{s \searrow 0} 
	\frac{
		\|I + s A_{ii}(t)\|_{\local i} - 1
	}{s} 
	|\xi_i(t, t_0, x)|_{\local i} \\
	&\quad\, + \sum_{j \neq i} \|A_{ij}(t)\|_{\local ij} |\xi_j(t, t_0, x)|_{\local j} \\
	&= \mu_{\local i}(A_{ii}(t)) |\xi_i(t, t_0, x)|_{\local i} 
	+ \sum_{j \neq i} \|A_{ij}(t)\|_{\local ij} |\xi_j(t, t_0, x)|_{\local j} \\
	&\leq \bar a^\lin_{ii}(t_1, t_0) |\xi_i(t, t_0, x)|_{\local i} 
	+ \sum_{j \neq i} \bar a^\lin_{ij}(t_1, t_0) |\xi_j(t, t_0, x)|_{\local j}.
\end{align*}
Therefore,
\[
	\rpdv{t}{\begin{bmatrix}
		|\xi_1(t, t_0, x)|_{\local 1} \\
		\vdots \\
		|\xi_m(t, t_0, x)|_{\local m}
	\end{bmatrix}} \leq \bar A^\lin_\Network(t_1, t_0) \begin{bmatrix}
		|\xi_1(t, t_0, x)|_{\local 1} \\
		\vdots \\
		|\xi_m(t, t_0, x)|_{\local m}
	\end{bmatrix}
\]
for almost all $ t \in [t_0, t_1) $.
Since the maps $ t \mapsto |\xi_i(t, t_0, x)|_{\local i} $ are absolutely continuous and $ \bar A^\lin_\Network(t_1, t_0) $ is a constant Metzler matrix over $ [t_0, t_1] $, the upper bound in \eqref{eq:inter-ltv-soln-upper} follows by the comparison principle for differential inequalities (cf. \cite[Prop.~1]{ArcakMaidens2018}).
\end{proof}

\begin{proof}[Proof of Lemma~\ref{lem:inter-soln-upper-reach}]
We use a variational representation of the difference between two trajectories (see, e.g., \cite[Sec.~2.5]{Vidyasagar2002}) to construct an interconnected LTV system and then apply Lemma~\ref{lem:inter-ltv-soln-upper}.

Fix arbitrary initial states $ x = (x_1, \ldots, x_m), \bar x = (\bar x_1, \ldots, \bar x_m) \in K $ with $ x_i, \bar x_i \in \R^{n_i} $, and define
\begin{equation*}
	\nu(t, \theta) 
	:= \theta \xi(t, t_0, \bar x) + (1 - \theta) \xi(t, t_0, x), 
	\qquad 
	t \geq t_0,\, \theta \in [0, 1].
\end{equation*}
Then
\begin{align*}
    &\quad\,\, \pdv{t}{(\xi(t, t_0, \bar x) - \xi(t, t_0, x))} \\
    &= f(t, \xi(t, t_0, \bar x)) - f(t, \xi(t, t_0, x)) \\
    &= f(t, \nu(t, 1)) - f(t, \nu(t, 0)) \\
    &= \int_{0}^{1} \jacob{x}{f(t, \nu(t, \theta))}\, \pdv{\theta}{\nu(t, \theta)} \dl{\theta} \\
    &= \left( \int_{0}^{1} \jacob{x}{f(t, \nu(t, \theta))} \dl{\theta} \right) 
    (\xi(t, t_0, \bar x) - \xi(t, t_0, x))
\end{align*}
for almost all $ t \geq t_0 $ 
(see, e.g., \cite[Sec.~4.2.4]{Liberzon2012} for a more detailed justification).
Thus, the trajectory difference $ \xi(t, t_0, \bar x) - \xi(t, t_0, x) $ is a solution of the interconnected LTV system \eqref{eq:inter-ltv} with initial state $ \bar x - x $ at $ t_0 $ and
\[
	A_{ij}(t) 
	:= \int_{0}^{1} \jacob{x_j}{f_i(t, \nu(t, \theta))} \dl{\theta} 
	\in \R^{n_i \times n_j}
\]
for $ i, j \in \{1, \ldots, m\} $.

Let $ \bar A^\lin_\Network(t, t_0) := \bar A_\Network^\reach(t, t_0, K) = [\bar a^\reach_{ij}(t, t_0, K)] $.
For each $ t \geq t_0 $, since the matrix measures $ \mu_{\local i}(\cdot) $ and induced norms $ \|\cdot\|_{\local ij} $ are convex functions, and $ \nu(t, \theta) \in \co(\xi(t, t_0, K)) $ for all $ \theta \in [0, 1] $, we have
\begin{align*}
	\sup_{s \in [t_0, t]} \mu_{\local i}(A_{ii}(s))
	&\leq \sup_{s \in [t_0, t]} 
	\int_{0}^{1} \mu_{\local i}(\jacob{x_i}{f_i(s, \nu(s, \theta))}) \dl{\theta} \\
	&\leq \sup_{s \in [t_0, t]} \max_{\theta \in [0, 1]} a_{ii}(s, \nu(s, \theta)) \\
	&\leq \bar a^\reach_{ii}(t, t_0, K), \\
	\sup_{s \in [t_0, t]} \|A_{ij}(s)\|_{\local ij}
	&\leq \sup_{s \in [t_0, t]} 
	\int_{0}^{1} \|\jacob{x_j}{f_i(s, \nu(s, \theta))}\|_{\local ij} \dl{\theta} \\
	&\leq \sup_{s \in [t_0, t]} \max_{\theta \in [0, 1]} a_{ij}(s, \nu(s, \theta))  \\
	&\leq \bar a^\reach_{ij}(t, t_0, K)
\end{align*}
for all diagonal and off-diagonal entries, respectively.
Hence \eqref{eq:inter-ltv-meas-max} holds.
Consequently, applying Lemma~\ref{lem:inter-ltv-soln-upper} over $ [t_0, t] $ yields \eqref{eq:inter-soln-upper-reach}.
\end{proof}

\subsection{Proof of Lemma~\ref{lem:grid-span}}\label{apx:grid-span}
\pushQED{\qed}
If \eqref{eq:grid-span-cond} holds, then
\begin{equation*}
    K \subset \bigcup_{x \in E(\theta)} (R(x, 2 \theta) \cap K) \subset \bigcup_{x \in E(\theta)} B_{f, t_0}(x, \varepsilon, T).
\end{equation*}
Thus, $ E(\theta) $ is a $ (T, \varepsilon) $-spanning set.
Consequently, the minimal cardinality of a $ (T, \varepsilon) $-spanning set satisfies
\begin{align*}
    S(f, t_0, \varepsilon, T, K) 
    &\leq \prod_{i=1}^{n} \left( 2 \left\lfloor \frac{r}{\theta_i} \right\rfloor + 1 \right) 
    \leq \prod_{i=1}^{n} \left( \frac{2 r}{\theta_i} + 1 \right).
\end{align*}
If \eqref{eq:grid-span-cond} holds for all $ T > 0 $ and $ \varepsilon > 0 $, then the definition of entropy \eqref{eq:ent-dfn} implies
\begin{align}
    h(f, t_0, K) 
    &\leq \limsup_{\varepsilon \searrow 0} 
    \limsup_{T \to \infty} 
    \sum_{i=1}^{n} \frac{\log(2 r/\theta_i + 1)}{T} \notag\\
    &\leq \limsup_{\varepsilon \searrow 0} 
    \limsup_{T \to \infty} 
    \sum_{i=1}^{n} \frac{\log(1/\theta_i)}{T} \notag\\
    &\quad\, + \sum_{i=1}^{n} \limsup_{\varepsilon \searrow 0} \limsup_{T \to \infty} \frac{\log(\theta_i + 2 r)}{T}, \label{eq:grid-span-ent-temp}
\end{align}
where the last inequality holds as the upper limit is a subadditive function. Moreover, each summand in the last term of \eqref{eq:grid-span-ent-temp} satisfies
\begin{align*}
    &\quad\, \limsup_{\varepsilon \searrow 0} \limsup_{T \to \infty} \frac{\log(\theta_i + 2 r)}{T} \\
    &\leq \limsup_{\varepsilon \searrow 0} \limsup_{T \to \infty} \frac{\max\{\log(2 \theta_i),\, \log(4 r)\}}{T} \\
    &= \max \left\{ \limsup_{\varepsilon \searrow 0} \limsup_{T \to \infty} \frac{\log\theta_i}{T},\, 0 \right\},
\end{align*}
where
the equality holds in part because $ r $ is independent of $ T $. Therefore, \eqref{eq:grid-span-ent-temp} implies \eqref{eq:grid-span-ent} under the condition \eqref{eq:grid-span-ent-cond}.
\popQED

\subsection{Proof of Lemma~\ref{lem:inter-ent-upper-eig-reach}}\label{apx:inter-ent-upper-eig-reach}
\pushQED{\qed}
Consider arbitrary initial states $ x = (x_1, \ldots, x_m), \bar x = (\bar x_1, \ldots, \bar x_m) \in K $ with $ x_i, \bar x_i \in \R^{n_i} $. 
By applying the subsystem-wise trajectory-separation bound in \eqref{eq:inter-soln-upper-reach}, together with the definition of the global norm $ |\cdot|_\Global $ in \eqref{eq:global-norm-dfn}, and using the monotonicity of the network norm $ |\cdot|_\Network $, we obtain
\begin{align*}
    &\quad\,\, |\xi(t, t_0, \bar x) - \xi(t, t_0, x)|_\Global \\
    &= \left| 
    	\begin{bmatrix}
			|\xi_1(t, t_0, \bar x) - \xi_1(t, t_0, x)|_{\local 1} \\
			\vdots \\
			|\xi_m(t, t_0, \bar x) - \xi_m(t, t_0, x)|_{\local m}
		\end{bmatrix} 
	\right|_\Network \\
	&\leq \left| 
		e^{\bar A_\Network^\reach(t, t_0, K) (t - t_0)} 
		\begin{bmatrix}
			|\bar x_1 - x_1|_{\local 1} \\
			\vdots \\
			|\bar x_m - x_m|_{\local m}
		\end{bmatrix} 
	\right|_\Network \\
	&\leq \|e^{\bar A_\Network^\reach(t, t_0, K) (t - t_0)}\|_\Network 
	\left| 
		\begin{bmatrix}
			|\bar x_1 - x_1|_{\local 1} \\
			\vdots \\
			|\bar x_m - x_m|_{\local m}
		\end{bmatrix} 
	\right|_\Network
\end{align*}
for all $ t \geq t_0 $. By the equivalence of norms on finite-dimensional vector spaces, there exist constants $ r_{\local 1}, \ldots, r_{\local m}, r_\Network > 0 $ such that
\begin{equation*}
	|v_i|_{\local i} \leq r_{\local i} |v_i|_\infty 
	\qquad 
	\forall v_i \in \R^{n_i},\, i \in \{1, \ldots, m\}
\end{equation*}
for the local norms $ |\cdot|_{\local i} $, and
\begin{equation*}
	|v|_\Network \leq r_\Network |v|_\infty 
	\qquad 
	\forall v \in \R^m
\end{equation*}
for the network norm $ |\cdot|_\Network $. 
Consequently, 
\begin{align*}
    &\quad\,\, |\xi(t, t_0, \bar x) - \xi(t, t_0, x)|_\Global \\
    &\leq \|e^{\bar A_\Network^\reach(t, t_0, K) (t - t_0)}\|_\Network 
    r_\Network 
    \max_{1 \leq i \leq m} |\bar x_i - x_i|_{\local i} \\
    &\leq \|e^{\bar A_\Network^\reach(t, t_0, K) (t - t_0)}\|_\Network 
    r_\Network 
    \max_{1 \leq i \leq m} r_{\local i} |\bar x_i - x_i|_\infty
\end{align*}
for all $ t \geq t_0 $. Therefore, for arbitrary $ T > 0 $ and $ \varepsilon > 0 $, we have
\begin{align*}
	&\quad\, \max_{t \in [t_0, t_0 + T]} 
	|\xi(t, t_0, \bar x) - \xi(t, t_0, x)|_\Global \\
	&\leq \left( 
		\sup_{t \in [t_0, t_0 + T]} 
		\|e^{\bar A_\Network^\reach(t, t_0, K) (t - t_0)}\|_\Network 
	\right) 
	r_\Network 
	\max_{1 \leq i \leq m} 
	r_{\local i} |\bar x_i - x_i|_\infty.
\end{align*}

Consider the set $ E(\theta) $ defined in Section~\ref{ssec:ent-span} with the vector $ \theta = (\bar\theta_1 \oneMx_{n_1}, \ldots, \bar\theta_m \oneMx_{n_m}) \in \R_{> 0}^n $, where
\[
    \bar\theta_i 
    := \varepsilon \left/ \left( 
    	2 r_\Network r_{\local i} 
    	\sup_{t \in [t_0, t_0 + T]} 
    	\|e^{\bar A_\Network^\reach(t, t_0, K) (t - t_0)}\|_\Network 
    \right) \right.,
\]
for $ i \in \{1, \ldots, m\} $, and $ \oneMx_{n_i} $ denotes the vector of ones in $ \R^{n_i} $.
The choice of $ \theta $ ensures that \eqref{eq:grid-span-cond} holds for 
the hyperrectangles $ R(x, 2\theta) $ defined in \eqref{eq:grid-rect}
and the open balls $ B_{f,t_0}(x,\varepsilon,T) $ defined in \eqref{eq:ball-dfn} with the global norm
$ |\cdot|_\Global $.
Hence Lemma~\ref{lem:grid-span} implies that $ E(\theta) $ is $ (T, \varepsilon) $-spanning. 
Moreover, since $ T > 0 $ and $ \varepsilon > 0 $ are arbitrary and all $ \bar\theta_i $ are nonincreasing in $ T $, the entropy bound in \eqref{eq:grid-span-ent} applies. 
Consequently,
\begin{align*}
    h(f, t_0, K) 
    &\leq \limsup_{\varepsilon \searrow 0} \limsup_{T \to \infty} 
    \sum_{i=1}^{m} \frac{n_i \log(1/\bar\theta_i)}{T} \\
    &= \limsup_{T \to \infty} \frac{n}{T} \sup_{t \in [t_0, t_0 + T]} 
    \log\|e^{\bar A_\Network^\reach(t, t_0, K) (t - t_0)}\|_\Network,
\end{align*}
where we use the fact that $ C_\varepsilon := \sum_{i=1}^{m} n_i \log(2 r_\Network r_{\local i}/\varepsilon) $ is independent of $ T $ and thus satisfies $ C_\varepsilon/T \to 0 $ as $ T \to \infty $. 
\popQED

\subsection{Proof of Lemma~\ref{lem:tv-vol-lower}}\label{apx:tv-vol-lower}
\pushQED{\qed}
For each $ v \in \R^n $, consider the LTV system \eqref{eq:ltv} with $ A(t) := \jacob{x}{f(t, \xi(t, t_0, v))} $.
By applying the standard variational equation for nonlinear time-varying systems (see, e.g., \cite[Th.~7, p.~24]{Coppel1965}) on intervals of continuity and joining the solutions at time discontinuities, 
its state-transition matrix is given by the solution Jacobian $ \jacob{x}{\xi(t, t_0, v)} $.

By Liouville's formula (see, e.g., \cite[Th.~4.1, p.~28]{Brockett1970}), the determinant satisfies
\[
    \det(\jacob{x}{\xi(t, t_0, v)}) 
    = e^{\int_{t_0}^{t} \tr(\jacob{x}{f(s, \xi(s, t_0, v))}) \dl{s}} 
    \qquad 
    \forall t \geq t_0.
\]
Therefore, by the change-of-variables formula,
\begin{align*}
    \vol(\xi(t, t_0, K)) 
    &= \int_{K} |{\det(\jacob{x}{\xi(t, t_0, v)})}| \dl{v} \\
    &\geq \left( \min_{v \in K} \det(\jacob{x}{\xi(t, t_0, v)}) \right) \vol(K) \\
    &= e^{\gamma(t, t_0)} \vol(K)
\end{align*}
for all $ t \geq t_0 $, which establishes \eqref{eq:tv-vol-lower}.
\popQED

\section{Proof sketches of Propositions~\ref{prop:tv-ent-upper-ini} and~\ref{prop:inter-ent-upper-eig-ini}}\label{apx:tv-inter-ent-upper-ini}
We first consider Proposition~\ref{prop:tv-ent-upper-ini}.
Its proof is analogous to that of Theorem~\ref{thm:tv-ent-bnd}.
The main difference lies in the construction of trajectory-separation bounds.
Specifically, we apply variational arguments to the line segment connecting two initial states, rather than the line segment connecting their corresponding trajectories at each time, as in the proof of Lemma~\ref{lem:inter-soln-upper-reach}.
Using Coppel's inequality (see, e.g., \cite[Th.~3, p.~47]{Vidyasagar2002}),
we obtain the following lemma.
The proof is given in Appendix~\ref{apx:tv-soln-upper-ini}.
\begin{lem}\label{lem:tv-soln-upper-ini}
For all initial states $ x, \bar x \in K $, the solutions of the nonlinear time-varying system \eqref{eq:tv} satisfy
\begin{equation}\label{eq:tv-soln-upper-ini}
     |\xi(t, t_0, \bar x) - \xi(t, t_0, x)| 
     \leq e^{\eta^\ini(t, t_0, K)} |\bar x - x| 
     \qquad 
     \forall t \geq t_0,
\end{equation}
where
\begin{equation}\label{eq:tv-jacob-meas-ini}
    \eta^\ini(t, t_0, K) 
    := \max_{v \in \co(K)} 
    \int_{t_0}^{t} \mu(\jacob{x}{f(s, \xi(s, t_0, v))}) \dl{s}.
\end{equation}
\end{lem}

Applying Lemma~\ref{lem:tv-soln-upper-ini} in place of Lemma~\ref{lem:inter-soln-upper-reach} leads to the constant $ \hat\mu^\iniinf $ defined in \eqref{eq:tv-meas-lim-ini} in Proposition~\ref{prop:tv-ent-upper-ini}, rather than the constant $ \hat\mu^\reachlim $ defined in \eqref{eq:tv-meas-lim-reach} in Theorem~\ref{thm:tv-ent-bnd}.
The complete proof of Proposition~\ref{prop:tv-ent-upper-ini} is given in Appendix~\ref{apx:tv-ent-upper-ini-full}.

We next consider Proposition~\ref{prop:inter-ent-upper-eig-ini}.
Its proof is analogous to that of Theorem~\ref{thm:inter-ent-upper-eig-reach}.
As in Proposition~\ref{prop:tv-ent-upper-ini}, we apply variational arguments to the line segment connecting two initial states, but use the subsystem-wise comparison argument of Lemma~\ref{lem:inter-ltv-soln-upper} instead of Coppel's inequality.
This leads to the following lemma.
The proof is given in Appendix~\ref{apx:inter-soln-upper-ini}.
\begin{lem}\label{lem:inter-soln-upper-ini}
Consider the interconnection matrix $ A_\Network(t, x) = [a_{ij}(t, x)] \in \R^{m \times m} $ defined in \eqref{eq:inter-sub-meas}.
Let $ \bar A_\Network^\ini(t, t_0, K) = [\bar a^\ini_{ij}(t, t_0, K)] \in \R^{m \times m} $ be any matrix-valued function satisfying
\begin{equation}\label{eq:inter-sub-meas-max-ini}
	\bar a^\ini_{ij}(t, t_0, K) 
	\geq \sup_{s \in [t_0, t]} 
	\max_{v \in \co(K)} 
	a_{ij}(s, \xi(s, t_0, v))
\end{equation}
for all $ t \geq t_0 $ and $ i, j \in \{1, \ldots, m\} $.
Then for all initial states $ x = (x_1, \ldots, x_m), \bar x = (\bar x_1, \ldots, \bar x_m) \in K $ with $ x_i, \bar x_i \in \R^{n_i} $, the solutions of the interconnected system \eqref{eq:inter} satisfy
\begin{align}
	&\begin{bmatrix}
		|\xi_1(t, t_0, \bar x) - \xi_1(t, t_0, x)|_{\local 1} \\
		\vdots \\
		|\xi_m(t, t_0, \bar x) - \xi_m(t, t_0, x)|_{\local m}
	\end{bmatrix} \notag\\
	&\qquad \leq e^{\bar A_\Network^\ini(t, t_0, K) (t - t_0)} 
	\begin{bmatrix}
		|\bar x_1 - x_1|_{\local 1} \\
		\vdots \\
		|\bar x_m - x_m|_{\local m}
	\end{bmatrix} 
	\qquad 
	\forall t \geq t_0. \label{eq:inter-soln-upper-ini}
\end{align}
\end{lem}

Applying Lemma~\ref{lem:inter-soln-upper-ini} in place of Lemma~\ref{lem:inter-soln-upper-reach} leads to the condition \eqref{eq:inter-sub-meas-sup-ini} for the comparison matrix $ \bar A^\iniinf_\Network $ in Proposition~\ref{prop:inter-ent-upper-eig-ini}, rather than the condition \eqref{eq:inter-sub-meas-lim-reach} for the comparison matrix $ \hat A^\reachlim_\Network $ in Theorem~\ref{thm:inter-ent-upper-eig-reach}.
The complete proof of Proposition~\ref{prop:inter-ent-upper-eig-ini} is given in Appendix~\ref{apx:inter-ent-upper-eig-ini-full}.

\section{Proof details of Propositions~\ref{prop:tv-ent-upper-ini} and~\ref{prop:inter-ent-upper-eig-ini}}\label{apx:tv-inter-ent-upper-ini-full}
\subsection{Proof of Lemma~\ref{lem:tv-soln-upper-ini}}\label{apx:tv-soln-upper-ini}
\pushQED{\qed}
For each $ v \in \R^n $, consider the LTV system \eqref{eq:ltv} with $ A(t) := \jacob{x}{f(t, \xi(t, t_0, v))} $.
By applying the standard variational equation for nonlinear time-varying systems (see, e.g., \cite[Th.~7, p.~24]{Coppel1965}) on intervals of continuity and joining the solutions at time discontinuities, 
its state-transition matrix is given by the solution Jacobian $ \jacob{x}{\xi(t, t_0, v)} $.

Fix arbitrary initial states $ x, \bar x \in K $, and define
\begin{equation*}
    \nu(t, \theta) 
    := \xi(t, t_0, \theta \bar x + (1 - \theta) x), 
    \qquad 
    t \geq t_0,\, 
    \theta \in [0, 1].
\end{equation*}
For each fixed $ \theta \in [0, 1] $, we have
\[
	\pdv{\theta}{\nu(t, \theta)} 
	= \jacob{x}{\xi(t, t_0, \theta \bar x + (1 - \theta) x)} 
	(\bar x - x).
\]
Thus, the derivative $ \pdv{\theta}{\nu(t, \theta)} $ is a solution of \eqref{eq:ltv} with initial state $ \bar x - x $ at $ t_0 $ and $ A(t) := \jacob{x}{f(t, \nu(t, \theta))} $. 

By Coppel's inequality (see, e.g., \cite[Th.~3, p.~47]{Vidyasagar2002}) and since $ \nu(t, \theta) \in \xi(t, t_0, \co(K)) $ for all $ \theta \in [0, 1] $, we have
\[
	|\pdv{\theta}{\nu(t, \theta)}| 
	\leq e^{\int_{t_0}^{t} \mu(\jacob{x}{f(s, \nu(s, \theta))}) \dl{s}} 
	|\bar x - x| 
	\leq e^{\eta^\ini(t, t_0, K)} 
	|\bar x - x|
\]
for all $ t \geq t_0 $.
Consequently,
\begin{align*}
    |\xi(t, t_0, \bar x) - \xi(t, t_0, x)| 
    &\leq \int_{0}^{1} |\pdv{\theta}{\nu(t, \theta)}| \dl{\theta} \\
    &\leq e^{\eta^\ini(t, t_0, K)} |\bar x - x|
\end{align*}
for all $ t \geq t_0 $, which establishes \eqref{eq:tv-soln-upper-ini}.
\popQED

\subsection{Proof of Proposition~\ref{prop:tv-ent-upper-ini}}\label{apx:tv-ent-upper-ini-full}
Combining the trajectory-separation bound in Lemma~\ref{lem:tv-soln-upper-ini} with the spanning set construction in Lemma~\ref{lem:grid-span} yields the following intermediate upper bound on entropy.

\begin{lem}\label{lem:tv-ent-upper-ini}
The topological entropy of \eqref{eq:tv} is upper bounded by
\begin{equation}\label{eq:tv-ent-upper-ini-temp}
    h(f, t_0, K) 
    \leq \limsup_{T \to \infty} \frac{n}{T} 
    \max_{t \in [t_0, t_0 + T]} \eta^\ini(t, t_0, K),
\end{equation}
where the function $ \eta^\ini(t, t_0, K) $ is defined in \eqref{eq:tv-jacob-meas-ini}.
\end{lem}

\begin{proof}
By the equivalence of norms on finite-dimensional vector spaces, there exists a constant $ r > 0 $ such that
\begin{equation*}
	|v| \leq r |v|_\infty \qquad \forall v \in \R^n
\end{equation*}
for the norm $ |\cdot| $. 
Consider arbitrary initial states $ x, \bar x \in K $. By applying the trajectory-separation bound in \eqref{eq:tv-soln-upper-ini}, we obtain
\begin{align*}
    |\xi(t, t_0, \bar x) - \xi(t, t_0, x)| 
    &\leq e^{\eta^\ini(t, t_0, K)} |\bar x - x| \\
    &\leq e^{\eta^\ini(t, t_0, K)} r |\bar x - x|_\infty
\end{align*}
for all $ t \geq t_0 $. Therefore, for arbitrary $ T > 0 $ and $ \varepsilon > 0 $, we have
\begin{align*}
	&\quad\, \max_{t \in [t_0, t_0 + T]} |\xi(t, t_0, \bar x) - \xi(t, t_0, x)| \\
	&\leq \max_{t \in [t_0, t_0 + T]} e^{\eta^\ini(t, t_0, K)} r |\bar x - x|_\infty.
\end{align*}

Consider the set $ E(\theta) $ defined in Section~\ref{ssec:ent-span} with the vector $ \theta = (\theta_1, \ldots, \theta_n) \in \R_{> 0}^n $, where
\begin{equation*}
    \theta_i 
    := \varepsilon \left/ \left( 
    	2 r \max_{t \in [t_0, t_0 + T]} e^{\eta^\ini(t, t_0, K)} 
    \right) \right., \qquad i \in \{1, \ldots, n\}.
\end{equation*}
The choice of $ \theta $ ensures that \eqref{eq:grid-span-cond} holds for 
the hyperrectangles $ R(x, 2\theta) $ defined in \eqref{eq:grid-rect}
and the open balls $ B_{f,t_0}(x,\varepsilon,T) $ defined in \eqref{eq:ball-dfn}.
Hence Lemma~\ref{lem:grid-span} implies that $ E(\theta) $ is $ (T, \varepsilon) $-spanning. 
Moreover, since $ T > 0 $ and $ \varepsilon > 0 $ are arbitrary and all $ \theta_i $ are nonincreasing in $ T $, the entropy bound in \eqref{eq:grid-span-ent} applies. 
Consequently,
\begin{align*}
    h(f, t_0, K) &\leq \limsup_{\varepsilon \searrow 0} 
    \limsup_{T \to \infty} 
    \sum_{i=1}^{n} \frac{\log(1/\theta_i)}{T} \\
    &= \limsup_{T \to \infty} \frac{n}{T} 
    \max_{t \in [t_0, t_0 + T]} \eta^\ini(t, t_0, K),
\end{align*}
where we use the fact that $ C_\varepsilon := n \log(2 r/\varepsilon) $ is independent of $ T $ and thus satisfies $ C_\varepsilon/T \to 0 $ as $ T \to \infty $. 
\end{proof}

We now establish the entropy bound in Proposition~\ref{prop:tv-ent-upper-ini}.
\begin{proof}[Proof of Proposition~\ref{prop:tv-ent-upper-ini}]
First, by the infimum in \eqref{eq:tv-meas-lim-ini}, for each $ \delta > 0 $, there exists $ t_1 \geq t_0 $ such that 
\begin{equation*}
	\limsup_{t \to \infty} 
	\max_{v \in \co(K_1)} 
	\mu(\jacob{x}{f(t, \xi(t, t_1, v))}) 
	\leq \hat\mu^\iniinf + \delta,
\end{equation*}
where $ K_1 := \xi(t_1, t_0, K) $.
Next, by the upper limit in \eqref{eq:tv-meas-lim-ini}, there exists a sufficiently large $ t_2 \geq t_1 $ such that
\begin{equation*}
	\max_{v \in \co(K_1)} 
	\mu(\jacob{x}{f(t, \xi(t, t_1, v))}) 
	\leq \hat\mu^\iniinf + 2 \delta 
	\qquad 
	\forall t \geq t_2.
\end{equation*}

Applying Lemma~\ref{lem:tv-ent-upper-ini} with the new initial time $ t_1 $ and initial set $ K_1 $ yields
\[
    h(f, t_1, K_1) 
    \leq \limsup_{T \to \infty} \frac{n}{T} 
    \max_{t \in [t_1, t_1 + T]} \eta^\ini(t, t_1, K_1).
\]
Moreover, for all $ T > t_2 - t_1 $, we have
\begin{align*}
    &\quad\, \max_{t \in [t_1, t_1 + T]} \eta^\ini(t, t_1, K_1) \\
    &\leq \max \left\{ 
    	C_{12},\,
    	\max_{t \in [t_2, t_1 + T]} \max_{v \in \co(K_1)} 
    	\int_{t_1}^{t} \mu(\jacob{x}{f(s, \xi(s, t_1, v))}) \dl{s}
    \right\} \\
    &\leq \max \left\{ 
    	C_{12},\,
    	C_{12}
    	+ \max_{t \in [t_2, t_1 + T]} 
    	(\hat\mu^\iniinf + 2 \delta) (t - t_2)
    \right\} \\
    &\leq C_{12} + \max\{\hat\mu^\iniinf + 2 \delta,\, 0\} (t_1 + T - t_2),
\end{align*}
where
\[
	C_{12} := \max_{t \in [t_1, t_2]} \max_{v \in \co(K_1)} 
    \int_{t_1}^{t} \mu(\jacob{x}{f(s, \xi(s, t_1, v))}) \dl{s}
\]
is independent of $ T $ and thus satisfies $ C_{12}/T \to 0 $ as $ T \to \infty $.
Hence
\[
	h(f, t_1, K_1) 
	\leq \limsup_{T \to \infty} \frac{t_1 + T - t_2}{T}\, 
	n \max\{\hat\mu^\iniinf + 2 \delta,\, 0\}.
\]
Therefore, by Lemma~\ref{lem:ini-time} and since $ (t_1 + T - t_2)/T \to 1 $ as $ T \to \infty $, we have
\[
	h(f, t_0, K) 
	\leq n \max\{\hat\mu^\iniinf + 2 \delta,\, 0\}.
\]
Since $ \delta > 0 $ is arbitrary, the upper bound in \eqref{eq:tv-ent-upper-ini} follows.
\end{proof}

\subsection{Proof of Lemma~\ref{lem:inter-soln-upper-ini}}\label{apx:inter-soln-upper-ini}
\pushQED{\qed}
For each $ v \in \R^n $, consider the LTV system \eqref{eq:ltv} with $ A(t) := \jacob{x}{f(t, \xi(t, t_0, v))} $.
By applying the standard variational equation for nonlinear time-varying systems (see, e.g., \cite[Th.~7, p.~24]{Coppel1965}) on intervals of continuity and joining the solutions at time discontinuities, 
its state-transition matrix is given by the solution Jacobian $ \jacob{x}{\xi(t, t_0, v)} $.

Fix arbitrary initial states $ x, \bar x \in K $, and define
\begin{equation*}
    \nu(t, \theta) 
    := \xi(t, t_0, \theta \bar x + (1 - \theta) x), 
    \qquad 
    t \geq t_0,\, 
    \theta \in [0, 1].
\end{equation*}
For each fixed $ \theta \in [0, 1] $, we have
\[
	\pdv{\theta}{\nu(t, \theta)} 
	= \jacob{x}{\xi(t, t_0, \theta \bar x + (1 - \theta) x)} 
	(\bar x - x).
\]
Thus, the derivative $ \pdv{\theta}{\nu(t, \theta)} $ is a solution of the interconnected LTV system \eqref{eq:inter-ltv} with initial state $ \bar x - x $ at $ t_0 $ and 
\[
	A_{ij}(t) := \jacob{x_j}{f_i(t, \nu(t, \theta))}, \qquad i, j \in \{1, \ldots, m\}.
\] 

Let $ \bar A^\lin_\Network(t, t_0) := \bar A_\Network^\ini(t, t_0, K) = [\bar a^\ini_{ij}(t, t_0, K)] $.
For each $ t \geq t_0 $, since $ \nu(t, \theta) \in \xi(t, t_0, \co(K)) $ for all $ \theta \in [0, 1] $, we have
\begin{align*}
	\sup_{s \in [t_0, t]} \mu_{\local i}(A_{ii}(s))
	&\leq \sup_{s \in [t_0, t]} \max_{v \in \co(K)} a_{ii}(s, \xi(s, t_0, v)) \\
	&\leq \bar a^\ini_{ii}(t, t_0, K), \\
	\sup_{s \in [t_0, t]} \|A_{ij}(s)\|_{\local ij}
	&\leq \sup_{s \in [t_0, t]} \max_{v \in \co(K)} a_{ij}(s, \xi(s, t_0, v))  \\
	&\leq \bar a^\ini_{ij}(t, t_0, K)
\end{align*}
for all diagonal and off-diagonal entries, respectively.
Hence \eqref{eq:inter-ltv-meas-max} holds.
Consequently, applying Lemma~\ref{lem:inter-ltv-soln-upper} over $ [t_0, t] $ gives 
\begin{align*}
	&\begin{bmatrix}
		|\pdv{\theta}{\xi_1(t, t_0, \theta \bar x + (1 - \theta) x)}|_{\local 1} \\
		\vdots \\
		|\pdv{\theta}{\xi_m(t, t_0, \theta \bar x + (1 - \theta) x)}|_{\local m}
	\end{bmatrix} \\
	&\qquad \leq e^{\bar A_\Network^\ini(t, t_0, K) (t - t_0)} 
	\begin{bmatrix}
		|\bar x_1 - x_1|_{\local 1} \\
		\vdots \\
		|\bar x_m - x_m|_{\local m}
	\end{bmatrix} 
	\qquad 
	\forall t \geq t_0.
\end{align*}
Integrating over $ \theta \in [0, 1] $ and applying the triangle inequality componentwise yields \eqref{eq:inter-soln-upper-ini}.
\popQED

\subsection{Proof of Proposition~\ref{prop:inter-ent-upper-eig-ini}}\label{apx:inter-ent-upper-eig-ini-full}
Combining the trajectory-separation bound in Lemma~\ref{lem:inter-soln-upper-ini} with the spanning set construction in Lemma~\ref{lem:grid-span} yields the following intermediate upper bound on entropy.
\begin{lem}\label{lem:inter-ent-upper-eig-ini}
Consider the interconnection matrix $ A_\Network(t, x) = [a_{ij}(t, x)] \in \R^{m \times m} $ defined in \eqref{eq:inter-sub-meas}.
Let $ \bar A_\Network^\ini(t, t_0, K) = [\bar a^\ini_{ij}(t, t_0, K)] \in \R^{m \times m} $ be any locally bounded matrix-valued function satisfying \eqref{eq:inter-sub-meas-max-ini} for all $ t \geq t_0 $ and $ i, j \in \{1, \ldots, m\} $.
Then the topological entropy of \eqref{eq:inter} is upper bounded by
\begin{equation}\label{eq:inter-ent-upper-eig-ini-temp}
    h(f, t_0, K) 
    \leq \limsup_{T \to \infty} \frac{n}{T} 
    \sup_{t \in [t_0, t_0 + T]} \log\|e^{\bar A^\ini_\Network(t, t_0, K) (t - t_0)}\|_\Network.
\end{equation}
\end{lem}

The proof of Lemma~\ref{lem:inter-ent-upper-eig-ini} follows that of Lemma~\ref{lem:inter-ent-upper-eig-reach},
with Lemma~\ref{lem:inter-soln-upper-ini} and $ \bar A_\Network^\ini(t, t_0, K) $ in place of
Lemma~\ref{lem:inter-soln-upper-reach} and $ \bar A_\Network^\reach(t, t_0, K) $, respectively.

We now establish the entropy bound in Proposition~\ref{prop:inter-ent-upper-eig-ini}. 
\begin{proof}[Proof of Proposition~\ref{prop:inter-ent-upper-eig-ini}]
By the infimum in \eqref{eq:inter-sub-meas-sup-ini}, for each $ \delta_1 > 0 $, there exists $ t_1 \geq t_0 $ such that
\begin{equation*}
	\sup_{t \geq t_1} \max_{v \in \co(\xi(t_1, t_0, K))} 
	a_{ij}(t, \xi(t, t_1, v)) 
	\leq \bar a^\iniinf_{ij} + \delta_1
\end{equation*}
for all $ i, j \in \{1, \ldots, m\} $.

Consider the interconnected system \eqref{eq:inter} with the new initial time $ t_1 $ and the corresponding initial set $ K_1 := \xi(t_1, t_0, K) $.
As $ \xi(t, t_0, K) = \xi(t, t_1, K_1) $, the matrix-valued function $ \bar A_\Network^\ini(t, t_1, K_1) := \bar A^\iniinf_\Network + \delta_1 \oneMx $, where $ \oneMx $ denotes the matrix of ones in $ \R^{m \times m} $, satisfies \eqref{eq:inter-sub-meas-max-ini} for all $ t \geq t_1 $.

The remainder of the proof follows that of Theorem~\ref{thm:inter-ent-upper-eig-reach}, with Lemma~\ref{lem:inter-ent-upper-eig-ini} and $ \bar A^\iniinf_\Network $ in place of Lemma~\ref{lem:inter-ent-upper-eig-reach} and $ \hat A^\reachlim_\Network $, respectively.
\end{proof}

\footnotesize
\bibliography{reference}

\end{document}